\documentclass[11pt]{article}

\usepackage{latexsym}
\usepackage{amsmath, amsfonts,amssymb, amsthm, euscript,makeidx,color,mathrsfs}
\usepackage{comment}

\newtheorem{assumption}{Assumption}
\def\qed{ \ \vrule width.2cm height.2cm depth0cm\smallskip}

\newcommand{\la}{\langle}
\newcommand{\ra}{\rangle}

\newcommand{\ba}{\begin{array}}
\newcommand{\ea}{\end{array}}
\newcommand{\be}{\begin{equation}}
\newcommand{\ee}{\end{equation}}
\newcommand{\bea}{\begin{eqnarray}}
\newcommand{\eea}{\end{eqnarray}}
\newcommand{\beaa}{\begin{eqnarray*}}
\newcommand{\eeaa}{\end{eqnarray*}}

\def\dbE{\mathbb{E}}
\def\dbF{\mathbb{F}}

\def\dbL{\mathbb{L}}

\def\dbP{\mathbb{P}}
\def\dbR{\mathbb{R}}
\def\dbS{\mathbb{S}}

\def\a{\alpha}
\def\b{\beta}
\def\g{\gamma}
\def\d{\delta}
\def\e{\varepsilon}

\def\l{\lambda}

\def\si{\sigma}

\def\f{\varphi}
\def\th{\theta}
\def\o{\omega}

\def\G{\Gamma}
\def\D{\Delta}
\def\Th{\Theta}
\def\L{\Lambda}

\def\O{\Omega}
\def\cA{{\cal A}}

\def\cF{{\cal F}}

\def\cL{{\cal L}}

\def\cN{{\cal N}}
\def\cO{{\cal O}}
\def\cP{{\cal P}}

\def\cU{{\cal U}}

\def\cZ{{\cal Z}}
\def\no{\noindent}

\def\bs{\bigskip}
\def\q{\quad}
\def\qq{\qquad}

\def\pa{\partial}
\def\cd{\cdot}
\def\cds{\cdots}

\def\td{\nabla}

\def\bx{{\bf x}}

\def\tr{\hbox{\rm tr}}

\def\qed{ \hfill \vrule width.25cm height.25cm depth0cm\smallskip}

\newcommand{\basa}{\begin{assumption}}
\newcommand{\easa}{\end{assumption}}

\newcommand{\bas}{\begin{assum}}
\newcommand{\eas}{\end{assum}}

\def\pa{\partial}

 \def\cd{\cdot}
\def\cds{\cdots}

\def\tr{\hbox{\rm tr$\,$}}

\def\ol{\overline}
\def\ul{\underline}

\def\dis{\displaystyle}

\def\bx{{\bf x}}

\def\1{{\bf 1}}

\def\:{\!:\!}
\def\reff#1{{\rm(\ref{#1})}}
\def \proof{{\noindent \bf Proof\quad}}

at 9pt

\newtheorem{thm}{Theorem}[section]
\newtheorem{lem}[thm]{Lemma}

\newtheorem{prop}[thm]{Proposition}
\newtheorem{rem}[thm]{Remark}
\newtheorem{eg}[thm]{Example}
\newtheorem{defn}[thm]{Definition}
\newtheorem{assum}[thm]{Assumption}

\begin{document}

\title{\bf  Forward Backward SDEs in Weak Formulation }

\author{
Haiyang Wang\thanks{ School of Mathematics and Statistics, Shandong Normal University, Jinan, 250014, China. Email: health.sea@163.com.  This author is supported in part by the Distinguished Middle-Aged and Young Scientist Encourage and Reward Foundation of Shandong Province (ZR2017BA033).} ~ and ~{Jianfeng Zhang}\thanks{\noindent
Department of Mathematics, University of Southern California, Los
Angeles, CA 90089, USA. E-mail: jianfenz@usc.edu. This author is
supported in part by NSF grant \#1413717. } \footnote{The second author would like to thank Daniel Lacker for very helpful discussion on Section \ref{sect-DPP}.}}

\date{\today}
\maketitle

\begin{abstract} Although having been developed for more than two decades, the theory of forward backward stochastic differential equations is still far from complete. In this paper, we take one step back and investigate the formulation of FBSDEs. Motivated from several considerations, both in theory and in applications, we propose to study FBSDEs in weak formulation, rather than the strong formulation in the standard literature. That is, the backward SDE is driven by the forward component, instead of by the Brownian motion noise. We establish the Feyman-Kac formula for FBSDEs in weak formulation, both in classical and in viscosity sense. Our new framework is efficient especially when the diffusion part of the forward equation involves the $Z$-component of the backward equation. 
\end{abstract}

\vfill \bs

\no{\bf Keywords.} \rm Forward backward SDEs, strong formulation, weak formulation, dynamic programming principle, stochastic maximum principle, quasilinear PDEs, path dependent PDEs, weak solution, viscosity solution, martingale problem.

\bs

\no{\it 2000 AMS Mathematics subject classification:} 60H07, 60H30, 35R60, 34F05

\eject

\section{Introduction}
\label{sect-Introduction}
\setcounter{equation}{0}

In the standard literature, a coupled FBSDE takes the following form:  
 \bea 
 \label{FBSDE-strong}
 \left\{\begin{array}{lll}
 \dis X_t= x+\int_0^tb(s,\Th_s)ds +\int_0^t \sigma(s,\Th_s)d B_s,\\
 \dis Y_t=g(X_T)+\int_t^Tf(s,\Th_s) ds -\int_t^TZ_sdB_s,
 \end{array}
 \right.\q t\in[0,T],\q \dbP_0\mbox{-a.s.} 
 \eea
where $\Th := (X, Y, Z)$ is the solution triplet, $B$ is a Brownian motion  under the probability measure $\dbP_0$ and the coefficients $b$, $\si$, $f$,  and $g$ are $\dbF^B$-progressively measurable in all variables. There have been many publications on the subject, see e.g. Antonelli \cite{Antonelli}, Ma, Protter \& Yong \cite{MPY},  Hu \& Peng \cite{HP},  Yong \cite{Yong1}, Peng \& Wu \cite{PW}, Pardoux \& Tang \cite{PTang}, Delarue \cite{Delarue},  Zhang \cite{Zhang2}, Ma, Wu, Zhang \& Zhang \cite{MWZZ}, as well as the monograph  Ma \& Yong \cite{MY-book}.  However, the theory is still far from complete. The existing methods in the literature provide quite different sets of sufficient conditions, and the unified approach proposed in \cite{MWZZ} works only in one dimensional case and the conditions there are rather technical. Even worse, many FBSDEs arising from applications do not fit in any existing works.

To understand the problem better, we take one step back and try to understand the formulation of the problem. Is \reff{FBSDE-strong} indeed the "right" formulation of the problem? As we will justify below, we feel the following alternative form, which we call FBSDEs in weak formulation, or simply weak FBSDEs, seems more appropriate in many situations:
 \bea 
 \label{FBSDE-weak}
 \left\{\begin{array}{lll}
 \dis X_t= x+\int_0^tb(s,\Th_s)ds +\int_0^t \sigma(s,\Th_s)d B_s;\\
 \dis Y_t=g(X_T)+\int_t^Tf(s,\Th_s) ds -\int_t^TZ_sdX_s,
 \end{array}
 \right. \q t\in[0,T],\q \dbP_0\mbox{-a.s.} 
 \eea
To indicate the difference, we denote  by $\Th^S:= (X^S, Y^S, Z^S)$ the solution to \reff{FBSDE-strong}  and $\Th^W:=(X^W, Y^W, Z^W)$ the solution to \reff{FBSDE-weak}, where the superscripts $^S$ and $^W$ stand for strong and weak, respectively. We note that in \reff{FBSDE-weak} the stochastic integration in the backward equation is against $dX_t$, not against the Brownian motion $dB_t$. In the case that 
\bea
\label{inverse}
\mbox{the mapping $z\mapsto z\si(t,x,y,z)$ has an inverse function $\psi(t,x,y,z)$,}
\eea 
by denoting $\tilde Z := Z^W \si(t, \Th^W_t)$ and thus $Z^W = \psi(t, X^W_t, Y^W_t, \tilde Z_t)$, one can easily check that $(X^W_t, Y^W_t, \tilde Z_t)$ is a solution to the following FBSDE in strong formulation:
 \bea 
 \label{FBSDE-tilde}
 \left\{\begin{array}{lll}
 \dis X_t= x+\int_0^t\tilde b(s,\Th_s)ds +\int_0^t \tilde \sigma(s,\Th_s)d B_s;\\
 \dis Y_t=g(X_T)+\int_t^T\tilde f(s,\Th_s) ds -\int_t^TZ_sdB_s,
 \end{array}
 \right.  \dbP_0\mbox{-a.s.} 
 \eea
 where, for $\th:=(t, x,y,z)$,
 \beaa
 \tilde b (\th) = b(t,x,y, \psi(\th)),\q   \tilde \si (\th) = \si(t,x,y, \psi(\th)),\q  \tilde f (\th) = f(t,x,y, \psi(\th)) - \psi(\th) \tilde b(\th).
 \eeaa
 When $\si = \si(t,x,y)$ is independent of $z$ and $\si>0$, it is clear that $\psi(t,\th)= {z\over \si(t,x,y)}$. However, when $\si$ depends on $z$, typically we do not have the inverse function $\psi$.
 
We justify the weak formulation \reff{FBSDE-weak} in four aspects. Firstly, in the option pricing and hedging theory, which is one of the main applications of BSDEs and FBSDEs, let $S$ denote the stock price driven by a Brownian motion $B$.  For a hedging portfolio $h$ with wealth value $V$, the self financing condition gives $d V_t = [\cds]dt + h_t dS_t$. Note that $(S, V)$ here correspond to $(X, Y)$ in FBSDE,  and the stochastic integration in $dV$ is against $dS_t$, not $dB_t$. In  simple models like Black-Scholes model, $S$ and $B$ generate the same filtration, then such difference is not crucial and there is no problem for using the strong formulation. However, for superhedging problem in incomplete markets, for example, one has to use $dS_t$ to superhedge, then the weak formulation is indeed more appropriate. In fact, in many practical applications, $X$ is the state process  we observe and $B$ is the noise used to model the distribution of $X$. Note that one rationale of using Brownian motion is the central limit theorem, where the convergence is in distribution, in this case  the value of $B$ may even not exist physically. So in these applications the weak formulation is more appropriate. 

Secondly, in Markovian setting and in the case $\si = \si(t,x,y)$, the FBSDE \reff{FBSDE-strong} is associated with the following  quasilinear PDE with terminal condition $u(T,x)=g(x)$:
\bea
\label{PDE-strong}
\pa_t u + {1\over 2} \si^2(t,x,u)\pa_{xx}^2 u + b(t,x,u, \pa_x u \si(t,x,u)) \pa_x u + f(t,x,u, \pa_x u\si(t,x,u))=0,
\eea
through the so called nonlinear Feynman-Kac formula:
\bea
\label{FK-strong}
Y^S_t = u(t, X^S_t),\q Z^S_t = \pa_x u(t, X^S_t) \si(t, X^S_t, u(t,X^S_t)).
\eea
However, when $\si$ depends on $Z$, the PDE will involve the inverse function $\psi$ in \reff{inverse} which typically does not exist. The weak formulation \reff{FBSDE-weak}, instead, corresponds to the following more natural PDE even in the case $\si = \si(t,x,y,z)$:
 \bea
\label{PDE-weak}
\pa_t u + {1\over 2} \si^2(t,x,u,\pa_x u) \pa_{xx}^2 u  + f(t,x,u, \pa_x u)=0,
\eea
and the nonlinear Feynman-Kac formula is also simpler:
\bea
\label{FK-weak}
Y^W_t = u(t, X^W_t),\q Z^W_t = \pa_x u(t, X^W_t).
\eea
In particular, in the option pricing and hedging theory, the representation \reff{FK-weak}  means exactly that $Z^W$ is the Delta-hedging portfolio $h$. 
The case that $\si$ depends on $Z$ indeed makes the difference between strong and weak formulations. For example, the following well known counterexample in strong formulation:
 \bea 
 \label{example-strong}
 X_t= x+\int_0^t Z_s d B_s;\q Y_t=X_T   -\int_t^TZ_sdB_s,
 \eea
has infinitely many solutions. However, the corresponding weak FBSDE is wellposed in the sense of Example \ref{eg-unique} and Remark \ref{rem-unique} below:
 \bea 
 \label{example-weak}
 X_t= x+\int_0^t Z_s d B_s;\q Y_t=X_T  -\int_t^TZ_sdX_s,
 \eea

Thirdly, as another major application, many FBSDEs arise from stochastic control problems through the stochastic maximum principle.  However, the stochastic control problem typically does not have optimal control in strong formulation. Indeed, even the following simple problem may not have an optimal control in strong formulation:
\bea
\label{control-strong}
X^\a_t = x + \int_0^t \a_s ds + B_t,\q V_0 := \sup_{\a\in \cU} \dbE^{\dbP_0}\Big[g(X^\a_\cd) + \int_0^T f(t, \a_t)dt\Big],
\eea
The corresponding control problem in weak formulation:
\bea
\label{control-weak}
\left.\ba{c}
X_t := x+ B_t,\q B^\a_t := B_t-\int_0^t \a_s ds,\q  d\dbP^\a := e^{\int_0^t \a_s dB_s -{1\over 2}\int_0^t |\a_s|^2ds} d\dbP_0,\\
 V_0 := \sup_{\a\in \cU} \dbE^{\dbP^\a}\Big[g(X_\cd) + \int_0^T f(t, \a_t)dt\Big],
\ea\right.
\eea
has optimal control under mild and natural conditions. Consequently, the associated FBSDE will have weak solution but no strong solution. It is more natural and convenient to write the FBSDE in weak formulation when one studies weak solutions.
 
 Fourthly, again for stochastic control problems, there are typically two approaches in the literature: the dynamic programming principle and the stochastic maximum principle. Both approaches lead to certain hamiltonians but the two hamiltonians for the same control problem look quite different.  As we observe, if one uses weak FBSDE as the adjoint equation involved in the stochastic maximum principle, then the hamiltonian will coincide with the one derived from the dynamic programming principle. In this sense, the weak formulation provides an intrinsic connection between the two approaches.

After carrying out the above motivations in details, we define weak solutions for weak FBSDEs and the equivalent forward backward martingale problems.  By utilizing the recently developed theory of path dependent PDEs, we establish the nonlinear Feynman-Kac formula for path dependent weak FBSDEs. That is, if the associate path dependent PDE has a classical solution, then the weak FBSDE has a (strong) solution.

Our main goal of this paper is to apply the viscosity solution method to establish the uniqueness of weak solution of  the weak FBSDE.   We shall follow the arguments in  Ma, Zhang, \& Zheng \cite{MZZ} and Ma \& Zhang \cite{MZ}, which study weak solutions for FBSDEs in strong formulation in the case that $\si$ is independent of $z$. Our arguments rely heavily on the regularity results for the PDE. Since such regularity results for the path dependent PDEs are not available in the literature, in this part we shall restrict to the Markovian case. The main idea is to study the so called nodal sets of the weak FBSDEs, whose upper and lower bounds provide viscosity subsolution and  supersolution of the PDE. Then, provided the comparison principle for viscosity solutions of the PDE, we obtain the uniqueness of weak solutions to the weak FBSDE. We remark that, as in \cite{MZZ, MZ},  the problem is equivalent to the so called martingale problem, see also Costantini \& Kurtz \cite{CK}  for the application of viscosity solution methods on martingale problems in an abstract framework. 

The rest of the paper is organized as follows. In Section \ref{sect-motivation} we motivate weak FBSDEs. In Section \ref{sect-classical} we define weak solutions and establish the nonlinear Feynman-Kac formula, provided the path dependent PDE has a classical solution. In Section \ref{sect-viscosity} we prove the existence and uniqueness of weak solutions for Markovian weak FBSDEs. Finally in Appendix we provide some counterexamples in control theory, which help to motivate the weak formulation, and provide some detailed arguments for the required regularities for the PDE.

\section{Some motivations for weak FBSDEs}
\label{sect-motivation}
\setcounter{equation}{0}
In this section we provide some heuristic motivations for weak FBSDE \reff{FBSDE-weak}. To simplify the presentation, we restrict to Markovian case in one dimensional setting. 

\subsection{Applications in option pricing and  hedging theory}
\label{sect-hedging}
Consider a financial market with a risky asset $S$ and a risk free asset with interest rate $r=0$ (for simplicity). Assume $S$ satisfies the following SDE:
\bea
\label{S}
S_t = S_0 + \int_0^t \si(s, S_s) dB_s,
\eea
where $B$ is a $\dbP_0$-Brownian motion (so we assume $\dbP_0$ is a risk neutral measure).  Given a portfolio $(\l, h)$ with value process $V_t = \l_t + h_t S_t$, the self-financing condition states that
\bea
\label{SF}
dV_t = h_t dS_t.
\eea
Now given an European type of option with payoff $\xi$ at terminal time $T$, we say a self-financing portfolio $(\l, h)$ is a hedging portfolio if $V_T = \xi$, $\dbP$-a.s. This, combining with \reff{SF}, leads to a backward SDE against $dS_t$:
\bea
\label{BSDEV}
V_t  = \xi - \int_t^T h_s dS_s.
\eea
Then  \reff{S}-\reff{BSDEV} become a decoupled weak FBSDE with solution $(X, Y, Z) = (S, V, h)$.
We remark that BSDE \reff{BSDEV}  can be rewritten in strong formulation:
\bea
\label{BSDEVstrong}
V_t = \xi - \int_t^T \tilde h_s dB_s,\q\mbox{where}\q \tilde h_t := h_t \si(t, S_t).
\eea
In particular, when $\si >0$, \reff{BSDEV} and \reff{BSDEVstrong} are equivalent.  This is why many papers in the literature could use the strong formulation.

The situation is different, however, in incomplete markets.  For example, consider the case that $S$ is scalar but $B$ is multi-dimensional.  Then $\si$ is a vector, and we assume that $\si$ is Lipschitz in $S$ so that \reff{S} has a strong solution $S$.  Assume further that we  observe the noise $B$ but can trade only $S$. Then $\xi$ is in general $\dbF^{B}$-measurable. By the martingale representation theorem, BSDE \reff{BSDEVstrong} always admits a solution $(V, \tilde h)$. However, since one cannot trade  $B$, the process $\tilde h$ is not a legitimate trading portfolio. For practical purpose  one has to solve the weak BSDE \reff{BSDEV}. In general $\tilde h$ may not be in the form of $h \si$, then in this case the strong BSDE \reff{BSDEVstrong} and the weak BSDE \reff{BSDEV} are not equivalent and in general the weak BSDE \reff{BSDEV} may not have a solution $(V, h)$.  One sensible resolution is to consider the super-hedging price:
\bea
\label{BSDE-superhedging}
V_0 := \inf\{y: \exists h ~\mbox{such that}~ y + \int_0^T h_s dS_s \ge \xi, ~\dbP_0\mbox{-a.s.}\}.
\eea
This is in the sprit of the weak FBSDE. Indeed, one can formulate it as a reflected BSDE in weak formulation, which is beyond the scope of this paper and is left for future research. 

An alternative explanation for the nonexistence of  solution to the weak FBSDE in  above situation is that $X$ does not have martingale representation property for $\dbF^B$-martingales. In this case,  for theoretical interest we may relax BSDE \reff{BSDEV} by applying the extended martingale representation theorem, see e.g. Protter \cite{Protter}:
\bea
\label{BSDEVN}
V_t = \xi - \int_t^T h_s dS_s + N_T-N_t,
\eea
where $N$ is an orthogonal martingale such that $N_0=0$ and $d\la S, N\ra_t = 0$. Then \reff{BSDEVN} will have a unique solution $(V, h, N)$.

\subsection{Nonlinear Feynman-Kac formula}
As is well known, in the case $\si = \si(t,x,y)$, the strong FBSDE \reff{FBSDE-strong} is associated with the quasilinear PDE \reff{PDE-strong} via the nonlinear Feynman-Kac formula \reff{FK-strong}. The problem becomes tricky when $\si = \si(t,x,y,z)$ because the PDE will involve the inverse function $\psi$ in \reff{inverse}, which typically does not exist. The weak FBSDE \reff{FBSDE-weak} is associated with the quasilinear PDE \reff{PDE-weak}, which is more natural at least in the following aspects:

$\bullet$ $\si$ may depend on $z$ and the PDE does not involve the inverse function $\psi$ in \reff{inverse}.

$\bullet$ The component $Z$ of the solution corresponds to $\pa_xu$ directly, rather than $\pa_x u ~ \!\!\si$. In particular, in the application to the option pricing and hedging theory, the $Z$ in weak formulation corresponds directly to the Delta-hedging portfolio.

$\bullet$ The PDE is more natural in the sense that the coefficients $\si, f$ depend directly on $\pa_x u$, instead of $\pa_x u \si$.

$\bullet$ It is more convenient to study weak solutions of the weak FBSDE, which is closely related to the viscosity solution of the PDE \reff{PDE-weak}, than that of the strong FBSDE.

To see the advantage of the weak formulation more directly in the case  that $\si$ depends on $z$,  let's consider the  counterexample \reff{example-strong}. It is well known that \reff{example-strong} has infinitely many solutions. Indeed, for any $Z \in \dbL^2(\dbF^B, \dbP_0)$, $X_t := Y_t := x+\int_0^t Z_s dB_s$ is a solution to \reff{example-strong}. However, the weak FBSDE \reff{example-weak} is wellposed in the following sense.
\begin{eg}
\label{eg-unique}
 The weak FBSDE \reff{example-weak} has a unique solution such that $Z\in \cZ := \{Z\in \dbL^4(\dbP_0): Z\neq 0\}$ .
 \end{eg}
Note that $Z\in \dbL^4(\dbP_0)$ implying $\dbE^{\dbP_0}[\int_0^T |Z_t|^2 d\la B\ra_t + \int_0^T |Z_t|^2 d\la X\ra_t] <\infty$, and thus $X, Y$ are $\dbP_0$-martingales. We shall comment on the requirement $Z\neq 0$ in Remark \ref{rem-unique} below.

\proof It is clear that 
\bea
\label{sol}
X_t = Y_t= x+ B_t,\qq Z_t = 1
\eea
is a solution to  \reff{example-weak}. We next show that it's  the unique solution such that $Z\in \cZ$.

For any  $(t,x, y)$ and $Z\in \dbL^4(\dbP_0)$, denote
\beaa
X^{t,x,Z}_s := x+ \int_t^s Z_r dB_r,\q Y^{t,x,y,Z}_s := y + \int_t^s Z_r dX^{t,x,Z}_r = y + \int_t^s |Z_r|^2 dB_r ,
\eeaa
and define 
\bea
\label{baru0}
\left.\ba{lll}
\ol u(t,x) := \inf\Big\{y: \exists Z\in \dbL^4(\dbP_0) ~\mbox{such that}~ Y^{t,x,y,Z}_T \ge X^{t,x,Z}_T,~\dbP_0\mbox{-a.s.}\Big\};\\
\ul u(t,x) := \sup\Big\{y: \exists Z\in \dbL^4(\dbP_0) ~\mbox{such that}~ Y^{t,x,y,Z}_T \le X^{t,x,Z}_T,~\dbP_0\mbox{-a.s.}\Big\}.
\ea\right.
\eea
Note that both $X^{t,x,Z}$ and $Y^{t,x,y,Z}$ are $\dbP_0$-martingales, then $Y^{t,x,y,Z}_T \ge X^{t,x,Z}_T$, $\dbP_0$-a.s. implies $y = \dbE^{\dbP_0}[Y^{t,x,y,Z}_T] \ge \dbE^{\dbP_0}[X^{t,x,Z}_T] = x$, and thus $\ol u(t,x)\ge x$. Similarly, $\ul u(t,x) \le x$ and thus $\ul u(t,x) \le x \le\ol u(t,x)$.  On the other hand, for any solution $(X, Y,Z)$ to \reff{FBSDE-weak}, by the definition of $\ol u(t,x)$ and $\ul u(t,x)$ we see that $\ol u(t,X_t) \le Y_t \le \ul u(t, X_t)$. Thus $\ol u(t,x) = \ul u(t,x) = u(t,x):= x$ and $Y_t = u(t, X_t)=X_t$. This implies further that $Z_t = |Z_t|^2$. Since $Z\neq 0$, we see that $Z=1$ and hence \reff{sol} is the unique solution.
\qed

 \begin{rem}
 \label{rem-unique}
 {\rm (i)  If we allow $Z=0$, then the solution is not unique. Indeed, for any $Z$ satisfying $Z=|Z|^2$ (namely $Z$ takes values $0$ and $1$), it is clear that $X_t = Y_t = x+\int_0^t Z_s dB_s$ is  a solution to weak FBSDE \reff{FBSDE-weak}. However, we note that even in this case, the relationship $Y_t = X_t$ still holds, and the decoupling function $u(t,x) = x$ is still unique. Moreover, without surprise, $u(t,x) = x$ is a solution to the PDE \reff{PDE-weak} corresponding to $b=0, \si = z, f=0$:
 \beaa
 \pa_t u + {1\over 2} |\pa_x u|^2 \pa^2_{xx} u = 0,\q u(T,x) = x.
 \eeaa
  
 (ii) When $Z=0$, this is exactly the case that $X$ has degenerate diffusion coefficient $\si$. As we will see in the paper, the nondegeneracy of $\si$ is crucial.
 
 (iii) As we mentioned in (i), even if we allow $Z=0$, the decoupling function $u(t,x)=x$ is still unique. However, when $X$ can be degenerate, $Y_t = u(t, X_t)$ and $d Y_t = Z_t dX_t$ do not imply $Z_t = \pa_x u(t, X_t) = 1$.  That's why the uniqueness fails in this degenerate case.
 \qed}
 \end{rem} 

To avoid the degeneracy issue, we may modify the example as follows.
\begin{eg}
\label{eg-unique2}
Let $\si>0$ be bounded such that the fixed point set $\cN := \{z: \si(z) = z\}\neq \emptyset$. 

(i) For any $\cN$-valued $Z\in \dbL^2(\dbF^B,\dbP_0)$, $Y_t := X_t := x+\int_0^t Z_s dB_s$ is a solution to the following strong FBSDE:
\beaa
X_t = x + \int_0^t \si(Z_s) dB_s,\q Y_t = X_T - \int_t^T Z_s dB_s.
\eeaa

(ii)  The corresponding weak FBSDE 
\beaa
X_t = x + \int_0^t \si(Z_s) dB_s,\q Y_t = X_T - \int_t^T Z_s dX_s.
\eeaa
has a unique solution 
\beaa
Y_t := X_t := x+ \si(1) B_t,\q Z_t := 1.
\eeaa
Here the uniqueness holds for $Z\in \dbL^2(\dbF^B,\dbP_0)$.
 \end{eg}
\proof (i) is obvious, and (ii) follows the same arguments as in Example \ref{eg-unique}.  In particular, the weak BSDE can be rewritten as:
\beaa
X_t = x + \int_0^t \si(Z_s) dB_s,\q Y_t = X_T - \int_t^T Z_s \si(Z_s) dB_s,
\eeaa
and then we see that $Z=1$ is the unique fixed point of: $\si(z)= z \si(z)$, thanks to the  nondegeneracy of $\si$. Moreover, since $\si$ is bounded, then $\dbE^{\dbP_0}[\int_0^T |Z_t|^2 d\la X\ra_t] <\infty$ for any $Z\in \dbL^2(\dbF^B,\dbP_0)$, so the uniqueness  holds for $Z\in \dbL^2(\dbF^B,\dbP_0)$.
\qed

\subsection{Connections with stochastic control theory}
\label{sect-control}

\subsubsection{Stochastic control in strong formulation}
\label{sect-SMP}
Consider a simple stochastic control problem in strong formulation:
\bea
\label{V0-strong}
&\dis V_0 := \sup_{\a\in \cA} V^\a_0,\q  \mbox{where}&\\
  &\dis  X^\a_t :=  \int_0^t b(s,\a_s) ds +\int_0^t \si(s,\a_s)d B_s,~ V^
\a_0:= \dbE^{\dbP_0}\Big[g(X^\a_T) + \int_0^T f(t, \a_t) dt\Big].&\nonumber
\eea
 Here the admissible controls $\a$ are $\dbF^B$-progressively measurable. Note that $V^\a_0 = Y^\a_0$, where
 \bea
 \label{controlYa}
 Y^\a_t = g(X^\a_T) + \int_t^T f(s, \a_s) ds - \int_t^T Z^\a_s dB_s.
 \eea
 
 We first use the stochastic maximum principle to derive an associated FBSDE. Let $\D \a$ be given such that $\a + \e \D \a\in \cA$ for any $\e\in [0, 1]$. Denote
 \beaa
 \td X^{\a,\D\a} := \lim_{\e\to 0} {1\over \e}[X^{\a+\e\D \a} - X^\a],\q \td V^{\a,\D \a}_0 :=  \lim_{\e\to 0} {1\over \e}[V^{\a+\e\D \a}_0 - V^\a_0].
 \eeaa
One can easily see that
 \beaa
 \td X^{\a,\D\a}_t &=& \int_0^t  b'(s,\a_s)\D \a_s ds + \int_0^t \si'(s,\a_s) \D \a_s dB_s,\\
  \td V^{\a,\D\a}_0 &=&\dbE^{\dbP_0}\Big[\pa_x g(X^\a_T) \td X^{\a,\D\a}_T + \int_0^T f'(t,\a_t) \D \a_tdt\Big],
 \eeaa
 where $b', \si', f'$ refer to  the derivatives with respect to  $\a$. 
 Introduce an adjoint BSDE:
 \bea
 \label{tildeYa}
 \tilde Y^\a_t = \pa_x g(X^\a_T) - \int_t^T \tilde Z^\a_s dB_s.
 \eea
 By applying It\^{o} formula on $\tilde Y^\a_t \td X^{\a,\D\a}_t$  we obtain
 \beaa
 \td V^{\a,\D\a}_0 = \dbE^{\dbP_0}\Big[ \int_0^T [\tilde Y^\a_t b'(t, \a_t) + \tilde Z^\a_t \si'(t, \a_t) +f'(t, \a_t)] \D \a_t dt\Big].
 \eeaa
 Now assume $\a^* \in \cA$ is an interior point of $\cA$ and is an optimal control.  Then $\td V^{\a^*, \D \a}_0 \le 0$ for arbitrary $\D \a$. This implies
 \bea
 \label{optimal}
\tilde Y^{\a^*}_t b'(t, \a^*_t) + \tilde Z^{\a^*}_t \si'(t, \a^*_t) +f'(t,  \a^*_t) = 0.
 \eea
 Assume further that \reff{optimal} determines an $\a^*$: $\a^*_t = I(t,  \tilde Y^{\a^*}_t, \tilde Z^{\a^*}_t)$ for a function $I$. Then  combining \reff{V0-strong}-\reff{tildeYa},  we obtain the following coupled FBSDE in strong formulation: 
  \bea
 \label{FBSDEcontrol}
 \left\{\ba{lll}
 \dis X_t = \int_0^t  b(s, I(s,  \tilde Y_s, \tilde Z_s)) ds + \int_0^t  \si(s, I(s,  \tilde Y_s, \tilde Z_s)) dB_s;\\
 \dis Y_t = g(X_T) + \int_t^T f(s,  I(s, \tilde Y_s, \tilde Z_s)) ds - \int_t^T Z_s dB_s;\\
 \dis \tilde Y_t = \pa_x g(X_T)  - \int_t^T \tilde Z_s dB_s.
 \ea\right.
 \eea
 However, the above FBSDE is typically not covered by the existing methods in the literature, especially since $\si$ depends on $\tilde Z$. 
 We remark that all the existing works on weak solutions of (strong) FBSDEs do not allow $\si$ depending on $Z$, see e.g.  Antonelli \& Ma \cite{AM}, Delarue \& Guatteri \cite{DG}, Ma, Zhang \& Zheng \cite{MZZ}, and Ma \& Zhang \cite{MZ}.
  
We thus turn to  weak FBSDE for which we can study weak solutions more conveniently. Rewrite the adjoint BSDE \reff{tildeYa} in the spirit of weak formulation:
 \bea
 \label{hatYa}
 \hat Y^\a_t = \pa_x g(X^\a_T) + \int_t^T  b(s,\a_s) \hat Z^\a_s ds - \int_t^T \hat Z^\a_s dX^\a_s.
 \eea
One can easily see that its solution is:  again assuming $\si>0$,
\bea
\label{hattildeYa}
\hat Y^\a_t:= \tilde Y^\a_t, \q  \hat Z^\a_t := \tilde Z^\a_t \si^{-1}(t,\a_t),
 \eea
and the optimality condition \reff{optimal} becomes
 \bea
 \label{optimal2}
\hat Y^{\a^*}_t b'(t, \a^*_t) + \hat Z^{\a^*}_t \si\si'(t, \a^*_t) +f'(t,\a^*_t) = 0.
 \eea
 Assume the above determines an optimal $\a^*$: $\a^*_t = \hat I(t, \hat Y^{\a^*}_t, \hat Z^{\a^*}_t)$ for a function $\hat I$. Then \reff{FBSDEcontrol} becomes a (multidimensional) FBSDE in weak formulation:
  \bea
 \label{FBSDEcontrol2}
 \left\{\ba{lll}
 \dis X_t = \int_0^t  b(s, \hat I(s, \hat Y_s, \hat Z_s)) ds + \int_0^t  \si(s, \hat I(s, \hat Y_s, \hat Z_s)) dB_s;\\
 \dis Y_t = g(X_T) + \int_t^T [f(s, \hat I(s, \hat Y_s, \hat Z_s)) +  b(s, \hat I(s, \hat Y_s, \hat Z_s)) Z_s] ds  - \int_t^T Z_s dX_s;\\
 \dis \hat Y_t = \pa_x g(X_T) + \int_t^T  b(s, \hat I(s, \hat Y_s, \hat Z_s)) \hat Z_sds - \int_t^T \hat Z_s dX_s.
 \ea\right.
 \eea
 
 \begin{rem}
\label{rem-optimal}
{\rm When the weak FBSDE \reff{FBSDEcontrol2} has no strong solution, but only weak solution, the stochastic optimization problem \reff{V0-strong} in strong formulation still does not have optimal control. To obtain the existence of optimal control, it is more appropriate to study the optimization problem in weak formulation, see Subsection \ref{ss-controlweak} below.
\qed}
\end{rem}

\subsubsection{Consistency with dynamic programming principle}
\label{sect-DPP}
 As is well known, another standard approach for stochastic control problem is the dynamic programming principle, which focuses more on the value function.  Assume the control $\a$ takes values in $A$. Then $V_0= u(0,0)$, where $u$ satisfies the following HJB equation:
 \bea
 \label{controlPDE}
 &\pa_t u + H(t, \pa_x u, \pa^2_{xx} u) =0,\q u(T,x) = g(x),&\\
 &\dis\mbox{where}~ H(t,z,\g) := \sup_{\a\in A} \Big[{1\over 2} \si^2(t,\a) \g + b(t,\a) z + f(t,\a)\Big].&\nonumber
 \eea

Assuming $u$ is sufficiently smooth and FBSDE \reff{FBSDEcontrol} is wellposed. By Yong \& Zhou \cite{YZ} Chapter 5, Theorem 4.1 we have    
\bea
\label{Yu1}
Y_t = u(t, X_t), ~ Z_t = \pa_x u(t, X_t) \si(t, \a^*_t),~ \tilde Y_t = \pa_x u(t, X_t), ~ \tilde Z = \pa^2_{xx} u(t, X_t) \si(t, \a^*_t),
\eea
 where $\a^*_t = I(t,\tilde Y_t, \tilde Z_t)$ is the optimal control. On the other hand, notice that the optimality condition \reff{optimal} can be viewed as the first order condition of
 \beaa
 \tilde H(t, \tilde y, \tilde z) := \sup_{\a\in A} \big[\tilde y b(t, \a) + \tilde z \si(t, \a) +f(t,\a) \big].
 \eeaa
However, we have the following  discrepancy which has already been noticed in  \cite{YZ}:
 \bea
 \label{Hdiff}
\tilde  H(t, X_t, \tilde Y_t, \tilde Z_t) &=& \pa_x u(t, X_t) b(t, \a^*_t)  + \pa^2_{xx} u(t, X_t) \si^2 (t, \a^*_t) + f(t, X_t, \a^*_t)\nonumber\\
 &\neq& H(t, X_t, \pa_x u(t, X_t) , \pa^2_{xx} u(t, X_t)).
  \eea

  This discrepancy is due to the fact that $\tilde Z$ involves $\si(t, \a)$ and thus twisted the optimization in the Hamiltonian.  It will disappear if we consider the weak FBSDE \reff{FBSDEcontrol2}. Indeed, in this case the optimality condition \reff{optimal2} can be viewed as the first order condition of
 \bea
 \label{Hhat}
 \hat H(t, \hat y, \hat z) := \sup_{\a\in A} \big[\hat y b(t, \a) + {1\over 2}\hat z \si^2(t, \a) +f(t,\a) \big].
 \eea
 Similar to \reff{Yu1} we have the correspondence for the solution to weak FBSDE \reff{FBSDEcontrol2}:
\bea
\label{Yu2}
Y_t = u(t, X_t), \q Z_t = \pa_x u(t, X_t),\q \hat Y_t = \pa_x u(t, X_t), \q \hat Z = \pa^2_{xx} u(t, X_t).
\eea
Then we have the desired identity:
  \bea
  \label{Hequal}
\hat H(t, \hat Y_t, \hat Z_t) &=& \pa_x u(t, X_t) b(t, \a^*_t)  + {1\over 2} \pa^2_{xx} u(t, X_t) \si^2 (t, \a^*_t) + f(t, \a^*_t)\nonumber\\
 &=& H(t, \pa_x u(t, X_t) , \pa^2_{xx} u(t, X_t)).
  \eea

\begin{rem}
\label{rem-I}
{\rm (i) It is clear that $\hat H = H$, with the correspondence $\hat y = z, \hat z = \g$. This is reflected in \reff{Yu2}.  In particular, we have $\hat Y_t = Z_t$ in this model.

(ii) The derivation of  \reff{optimal2} requires the differentiation of the coefficients $b, \si, f$ in $\a$. However, such differentiation is not needed for the optimization of  the Hamiltonian in  \reff{Hhat}. In fact, one may determine  $\hat I$ by the optimal arguments in  \reff{Hhat}, and then formally derive the same FBSDE \reff{FBSDEcontrol2}. These arguments are in the line of dynamic programming principle, rather than stochastic maximum principle.
\qed}
\end{rem}

   \subsubsection{Stochastic drift control under weak formulation}
   \label{ss-controlweak}
   
 To understand the weak FBSDE \reff{FBSDEcontrol2} better, we consider a special case that 
 \beaa
 \si = 1.
 \eeaa
 The general case with diffusion control will involve the second order BSDE introduced in Soner, Touzi, \& Zhang \cite{STZ}. 
 In this case \reff{optimal2} becomes
 \bea
 \label{optimal3}
\hat Y^{\a^*}_t b'(t, \a^*_t) +f'(t,\a^*_t) = 0.
 \eea
 Then the optimal control takes the form $\a^*_t = \hat I(t, \hat Y^{\a^*}_t)$ and thus  FBSDE \reff{FBSDEcontrol2} becomes
  \bea
 \label{FBSDEcontrol3}
 \left\{\ba{lll}
 \dis X_t = \int_0^t  b(s, \hat I(s, \hat Y_s)) ds +B_t;\\
 \dis Y_t = g(X_T) + \int_t^T [f(s, \hat I(s, \hat Y_s)) +  b(s, \hat I(s, \hat Y_s)) Z_s] ds  - \int_t^T Z_s dX_s;\\
 \dis \hat Y_t = \pa_x g(X_T) + \int_t^T  b(s, \hat I(s, \hat Y_s)) \hat Z_sds - \int_t^T \hat Z_s dX_s.
 \ea\right.
 \eea
Recalling \reff{Yu2} that $\hat Y = Z$, the second equation in \reff{FBSDEcontrol3} is equivalent to
\bea
\label{FBSDEcontrol4}
Y_t = g(X_T) + \int_t^T [f(s, \hat I(s, Z_s)) +  b(s, \hat I(s, Z_s)) Z_s] ds  - \int_t^T Z_s dX_s.
\eea
Moreover, note that \reff{optimal3} is the first order condition of the following optimization problem:
\bea
\label{f*}
f^*(t,z) := \sup_{\a\in A} [z b(t,\a) + f(t,\a)].
\eea
 Then, together with \reff{FBSDEcontrol3} and under appropriate technical conditions, \reff{FBSDEcontrol4}  leads to
\bea
 \label{FBSDEcontrol5}
 \left\{\ba{lll}
 \dis X_t = \int_0^t  I(s, Z_s) ds + B_t;\\
 \dis Y_t = g(X_T) + \int_t^T  f^*(s,  Z_s) ds - \int_t^T Z_s dX_s.
  \ea\right.
 \eea

The FBSDE \reff{FBSDEcontrol5} can be understood a lot easier  if we use weak formulation for the control problem:
 \bea
 \label{V0-weak}
 &\dis \bar V_0 := \sup_{\a\in \cA}  \bar V^\a_0 :=  \sup_{\a\in \cA}  \dbE^{\dbP^\a}\Big[g(X_T) + \int_0^T f(t,\a_t) dt\Big]&\\
  &\dis  \mbox{where}\q X_t := B_t,\q d\dbP^\a := \exp\Big(\int_0^T \a_t dB_t - {1\over 2} \int_0^T|\a_s|^2 dt\Big)d\dbP_0.&\nonumber
 \eea
 Note that $\bar V^\a_0 = \bar Y^\a_0$, where,  $B^\a_t := B_t - \int_0^t \a_s ds$ is a $\dbP^\a$-Brownian motion and
  \bea
 \label{controlYa2}
 \bar Y^\a_t &=& g(X_T) + \int_t^T f(s,\a_s) ds - \int_t^T \bar Z^\a_s dB^\a_s \nonumber\\
 &=& g(X_T) + \int_t^T [f(s,\a_s) + \a_s\bar Z^\a_s] ds - \int_t^T \bar Z^\a_s dX_s.
 \eea
 Consider the BSDE
 \bea
 \label{controlBSDE}
 \bar Y_t =  g(X_T) + \int_t^T f^*(s, \bar Z_s) ds - \int_t^T \bar Z_s dX_s,\q \dbP_0\mbox{-a.s.}
 \eea
By comparison of BSDE we see immediately that $\bar V_0 = \bar Y_0$ and $\a^*_t := \hat I(t, \bar Z_t)$ is an optimal control of \reff{V0-weak}, for the same $\hat I$ in \reff{FBSDEcontrol5}. Now together with the definition of $B^\a$ and $X = B$, we may rewrite \reff{controlBSDE} as 
\bea
 \label{XbarY}
 \left\{\ba{lll}
 \dis X_t = \int_0^t \hat I(s, \bar Z_s) ds + B^{\a^*}_t,\\
 \dis \bar Y_t =  g(X_T) + \int_t^Tf^*(s, \bar Z_s) ds - \int_t^T \bar Z_s dX_s.
 \ea\right.
 \eea
 In the spirit of weak solution as we will introduce in the next section, this is equivalent to \reff{FBSDEcontrol5}. So in this sense, the weak FBSDE \reff{FBSDEcontrol5}, or the more general one \reff{FBSDEcontrol2}, is more in the spirit of weak formulation. 

\begin{rem}
\label{rem-control}
{\rm (i) Recall \reff{V0-strong} with $\si=1$ and \reff{V0-weak}. Note that formally  $(B^\a, B, \dbP^\a)$ in weak formulation corresponds to $(B, X^\a, \dbP_0)$ in strong formulation. However, for fixed $\a$, note that in general $\a$ has different distributions under $\dbP^\a$ and under $\dbP_0$, so $V^\a_0 \neq \bar V^\a_0$.  But nevertheless, under appropriate conditions, their optimal values are equal: $V_0 = \bar V_0$. See more discussions along this line in Zhang \cite{Zhang-book} Chapter 9.

(ii) There are many situations that the optimal control in weak formulation exists but that in strong formulation does not. See some examples in Appendix.

(iii) The difference between strong formulation and weak formulation becomes more crucial when one considers zero sum stochastic differential games, see Hamadene \& Lepeltier \cite{HL} and Pham \& Zhang \cite{PZ}.
\qed}
\end{rem}

 \section{Weak solutions of FBSDEs and Feynman-Kac formula}
\label{sect-classical}
\setcounter{equation}{0}

Our objective is the following  weak FBSDE: 
 \bea
 \label{FBSDE}
 \left\{\ba{lll}
 \dis X_t = x + \int_0^t b(s, X_\cd, Y_s, Z_s) ds + \int_0^t \si(s, X_\cd, Y_s, Z_s) dB_s,\\
 \dis Y_t = g(X_\cd) + \int_t^T f(s, X_\cd, Y_s, Z_s) ds - \int_t^T Z_s dX_s + N_T-N_t,
 \ea\right. \dbP_0\mbox{-a.s.}
 \eea
 Here $(B, X, Y)$ take values in $\dbR^{d_0}\times \dbR^{d_1}\times \dbR^{d_2}$, and all other processes and functions have appropriate dimensions.  The coefficients $b, \si, f, g$ may depend on the paths of $X$, among them $b, \si, f$ are $\dbF^X$-progressively measurable in all variables, and $g$ is $\cF^X_T$-measurable. 
 
 Given a probability space $(\O, \dbF, \dbP)$, let  $\dbL^0(\dbF)$ denote the set of $\dbF$-progressively measurable processes with appropriate dimensions. For $p, q \ge 1$, denote
 \beaa
&\dis \dbL^{p, q} (\dbF, \dbP) := \big\{X\in \dbL^0(\dbF):  \dbE^\dbP\Big[ \big(\int_0^T |X_t|^p dt\big)^{q\over p}\big] <\infty\Big\},\q \dbL^p(\dbF, \dbP) := \dbL^{p,p}(\dbF, \dbP); &\\
&\dis \dbS^p(\dbF, \dbP) :=\big\{X \in \dbL^0(\dbF):  \mbox{$X$ is continuous, $\dbP$-a.s. and}~ \dbE^\dbP\big[ \sup_{0\le t\le T} |X_t|^p\big] < \infty\big\}.&
\eeaa
 Throughout this paper, we shall assume 
 \begin{assum}
 \label{assum-basic} (i)  $b, \si$ are bounded.
 
 (ii) $ f(t, \bx, 0, 0)$, $g(\bx)$ have polynomial growth in $\|\bx\| := \sup_{0\le t\le T}|\bx_t|$, and $f$ is uniformly Lipschitz continuous in $(y,z)$. 
 
 (iii)  $\si \si^\top \ge c^2_0 I_{d_1}$ as $d_1 \times d_1$-matrice, for some constant $c_0>0$.
 \end{assum}

 \begin{rem}
 \label{rem-path}
 {\rm   (i) For strong FBSDEs, typically the coefficients may depend on $B$. For weak FBSDEs, both for practical considerations and for theoretical reasons, it is more natural that the coefficients depend on $X$.  However, in a more general setting, for example in the incomplete market with observable noise as in Subsection \ref{sect-hedging}, we may allow the coefficients to depend on both $X$ and $B$. The problem will become harder in this case. In this paper we restrict to the case that the coefficients do not depend on $B$.
 
 (ii) As explained in Section \ref{sect-hedging}, the presence of $N$ is due to the fact that $X$ may not satisfy the martingale representation property. 
 \qed}
\end{rem}

 \subsection{Definitions}
 \label{sect-defn}
 
We introduce the following types of solutions. Recall $\Th = (X, Y, Z)$.
\begin{defn}
\label{defn-weak}
(i) We say a filtered probability space $(\O, \dbF, \dbP)$ and a quintuple of $\dbF$-progressively measurable processes $(B, \Th, N)$ is a weak solution of the weak FBSDE \reff{FBSDE} if $B$ is a $\dbP$-Brownian motion, $N\in \dbS^2(\dbF, \dbP)$ is a $\dbP$-martingale orthogonal to $X$ with $N_0=0$, $X, Y \in \dbS^2(\dbF, \dbP)$, $Z\in \dbL^{2}(\dbF,\dbP)$, and \reff{FBSDE} holds $\dbP$-a.s. 

(ii) We say a weak solution is semi-strong if $(Y, Z)$ are $\dbF^X$-progressively measurable.

(iii) We say a weak solution is strong if $N=0$ and $\Th$ is $\dbF^B$-progressively measurable.
\end{defn}

Given our conditions, all weak solutions actually have stronger integrability.
\begin{lem}
\label{lem-integrability}
Let Assumption \ref{assum-basic} hold and $(B, \Th, N, \dbP)$ be a weak solution to \reff{FBSDE}. Then
\beaa
 \dbE^\dbP\Big[\sup_{0\le t\le T} \big[|X_t|^p + |Y_t|^p + |N_t|^p\big] + \big(\int_0^T |Z_t|^2dt\big)^{p\over 2}\Big] < \infty,~\mbox{for any}~p\ge 1.
 \eeaa
 \end{lem}
\proof By the boundedness of $b, \si$, the estimate for $X$ is obvious. Since $f(t, \bx, 0, 0)$ and $g(\bx)$ have polynomial growth, we have $\dbE^\dbP\Big[|g(X_\cd)|^p + \int_0^T |f(t, X_\cd)|^pdt\Big] <\infty$. Now by the uniform Lipschitz continuity of $f$ in $(y,z)$, the rest estimates follows from standard BSDE arguments, see e.g.  El Karoui \& Huang \cite{KH}.  
\qed

As in Stroock \& Varadahn \cite{SV}, weak solutions are closely related to martingale problems. Motivated by  Ma, Zhang \& Zheng \cite{MZZ} which studies strong FBSDE with $\si$ independent of $z$, we introduce the following forward-backward martingale problem.
  
\begin{defn}
\label{defn-mg}
{\rm Let $\O := C([0, T], \dbR^{d_1}) \times C([0, T], \dbR^{d_2})$   be the canonical space,  $(X, Y)$ the canonical processes, and $\dbF=\dbF^{X, Y}$ the natural filtration.   We say $(\dbP, Z)$ is a solution to the forward-backward martingale problem of \reff{FBSDE} if:

(i) $\dbP(X_0=x)=\dbP(Y_T = g(X_\cd)) = 1$ and $X, Y \in \dbS^2(\dbF, \dbP)$, $Z\in \dbL^2(\dbF,\dbP)$.

(ii) The following two processes are $\dbP$-martingales:
\beaa
& \dis M^X_t := X_t - \int_0^t b(s, X_\cd, Y_s, Z_s)ds,&\\
 &\dis M^Y_t := Y_t + \int_0^t f(s, X_\cd, Y_s, Z_s)ds - \int_0^t Z_s b(s, X_\cd, Y_s, Z_s)ds.
\eeaa

(iii) $d\la M^X\ra_t = \si\si^\top (t, X_\cd, Y_t, Z_t) dt$ and  $d \la M^Y, M^X\ra_t = Z_t d\la X\ra_t$, $\dbP$-a.s.
\qed}
\end{defn}

\begin{prop}
\label{prop-equiv}
Let Assumption \ref{assum-basic} hold. Then a weak solution to FBSDE \reff{FBSDE} is equivalent to a solution to the forward-backward martingale problem of \reff{FBSDE}.
\end{prop}
\proof Let $(\O, \dbF, \dbP, B, \Th, N)$ be a weak solution to FBSDE \reff{FBSDE}. Note that $d \la Y, X\ra_t = Z_t d\la X\ra_t = Z_t \si \si^\top(t, X_\cd, Y_t, Z_t)dt$ and $\la X\ra, \la Y\ra$ are all $\dbF^{X, Y}$-progressively measurable.  Since $\si\si^\top >0$, then $Z$ is also $\dbF^{X, Y}$-progressively measurable. Now by recasting everything into the canonical space of $(X, Y)$,  it is straightforward to verify that $(\dbP, Z)$  is a solution to the forward-backward martingale problem of \reff{FBSDE}.

To see the other direction, let $(\O, \dbF, X, Y)$ be the canonical  setting in Definition \ref{defn-mg} and $(\dbP, Z)$  a solution to the forward-backward martingale problem of \reff{FBSDE}.  Note that Assumption  \ref{assum-basic} (iii) implies $d_0 \ge d_1$, and there exist orthogonal matrices $U\in \dbR^{d_1\times d_1}$ and $V\in \dbR^{d_0\times d_0}$ as well as $k_1, \cds, k_{d_1} \neq 0$ such that 
\beaa
\si (t, X_\cd, Y_t, Z_t)= U_t~ [ K_t, 0] ~V_t\q\mbox{where $K$ is the diagonal matrix of $k_1, \cds, k_{d_1}$},
\eeaa
and $0$ refers to the $d_1 \times (d_0-d_1)$-zero matrix. It is clear that $U, V, K$ are $\dbF$-progressively measurable processes. Denote
\beaa
\tilde B_t := \int_0^t K_s^{-1} U_s^\top dM^X_s.
\eeaa
Then $\tilde B$ is a continuous local martingale under $\dbP$ and
\beaa
d \la \tilde B \ra_t =  K_t^{-1} U_t^\top \si \si^\top U_t K_t^{-1} dt =  K_t^{-1} U_t^\top U_t [ K_t, 0] V V^\top [ K_t, 0]^\top U_t^\top U_t  K_t^{-1} dt = I_{d_1} dt.
\eeaa
By Levy's characterization theorem we see that $\tilde B$ is a $\dbP$-Brownian motion. Now let $\bar B$ be an $d_0-d_1$-dimensional Brownian motion independent of $\dbF$, and let us extend $\dbF$ to $\hat \dbF := \dbF \vee \dbF^{\bar B}$, and still denote the probability measure as $\dbP$. Then $\hat B := [B^\top, \bar B^\top]^\top$ is a $d_0$-dimensional $\dbP$-Brownian motion. Thus
\beaa
&\dis dM^X_t = U_t  K_t d \tilde B_t = U_t [ K_t, 0] d\hat B_t= \si(t, X_\cd, Y_t, Z_t) d B_t,
\eeaa
where $d B_t := V_t^\top d \hat B_t$ is also a $d_0$-dimensional $\dbP$-Brownian motion, since $V$ is orthogonal. 
Now define
\bea
\label{BN}
 N_t := Y_0-M^Y_t  + \int_0^t Z_s dM^X_s.
\eea
Then  $N$ is a $\dbP$-martingale. Note that $d\la X, N\ra_t = - d\la X, M^Y\ra_t + Z_t d\la X\ra_t=0$. Then $(\O, \hat\dbF, \dbP, B, X, Y, Z, N)$ is a weak solution to FBSDE \reff{FBSDE}. 
\qed

\begin{rem}
\label{rem-sigma}
{\rm Note that the martingale problem involves only $\si\si^\top$, not the $\si$ itself. Then by Proposition \ref{prop-equiv} we may assume without loss of generality that
\bea
\label{assum-si}
d_0 = d_1 =: d,\q \si ~\mbox{is symmetric and}~ \si \ge c_0I_d.
\eea
In the rest of the paper this will be enforced.
\qed}
\end{rem}

Given \reff{assum-si}, we  have another equivalence result.

\begin{prop}
\label{prop-b}
Let Assumption \ref{assum-basic} hold. Then  FBSDE \reff{FBSDE} admits a weak solution if and only if \reff{FBSDE} with coefficients $(0, \si, f, g)$ has a weak solution.
\end{prop}
\proof We assume without loss of generality that \reff{assum-si} holds. Let $(B, \Th, N, \dbP)$ be a weak solution to FBSDE \reff{FBSDE} with coefficients $(b, \si, f, g)$. Denote 
\beaa
\th_t := -\si^{-1} b(t, X_\cd, Y_t, Z_t),~ \tilde B_t := B_t - \int_0^t \th_s ds,~ d\tilde \dbP := \exp(\int_0^t \th_s dB_s - {1\over 2} \int_0^t |\th_s|^2 ds) d\dbP.
\eeaa
Then $\th$ is bounded, and thus it follows from Lemma \ref{lem-integrability}  that $(\Th, N)$ have the desired integrability under $\tilde \dbP$. 
Since $B$ and $N$ are orthogonal, by Girsanov theorem one can easily check that $(\tilde B, \Th, N, \tilde \dbP)$ is a weak solution to FBSDE \reff{FBSDE} with coefficients $(0, \si, f, g)$. This proves the only if part. The if part can be proved similarly.
\qed

\subsection{Path dependent PDEs}
In this subsection we introduce the PPDE in the setting of Ekren, Touzi, \& Zhang \cite{ETZ1, ETZ2}.  Let $\O := C([0, T], \dbR^{d})$   be the canonical space equipped with $\dis\|\o\|:= \sup_{0\le t\le T} |\o_t|$,  $X$ the canonical process, $\dbF:=\dbF^X$ the natural filtration, and $\L:= [0, T]\times \O$  equipped with 
\beaa
{\bf d}((t,\o), (t',\o')):= |t-t'|+ \sup_{0\le s\le T} |\o_{t\wedge s} - \o'_{t'\wedge s}|.
\eeaa
For some generic dimension $m$,  let $C^0(\L; \dbR^m)$ be the space of continuous functions $\L\to \dbR^m$.

Next, let $\cP$ denote the set of semimartingale measures $\dbP$ whose drift and diffusion characteristics are bounded, 
and $C^{1,2}(\L; \dbR)$ be the  space of $u \in  C^0(\L; \dbR)$ such that there exist $\pa_t u \in C^0(\L; \dbR)$, $\pa_\o u \in C^0(\L; \dbR^{1\times d})$ (row vector for convenience!), and symmetric $\pa^2_{\o\o} u \in  C^0(\L; \dbR^{d \times d})$ satisfying: for all $\dbP\in \cP$, $u(t, X_\cd)$ is a semimartingale and the following {\it functional It\^{o} formula} holds:
\bea
\label{FIto}
d u(t, X_\cd) = \pa_t u(t, X_\cd) dt + \pa_\o u(t, X_\cd)  dX_t + {1\over 2}\pa^2_{\o\o} u(t, X_\cd) : d \la X\ra_t,\q \dbP\mbox{-a.s.}
\eea
For each $u\in C^{1,2}(\Th; \dbR)$,  by \cite{ETZ1} the path derivatives $\pa_t u, \pa_\o u, \pa^2_{\o\o} u$ are unique.  Moreover, we say $u = [u_1,\cds, u_m]^\top \in C^{1,2}(\L; \dbR^m)$ if each $u_i  \in C^{1,2}(\L; \dbR)$ for $i=1,\cds, m$.

Denote $f = [f_1,\cds, f_{d_2}]^\top$. The weak FBSDE \reff{FBSDE} is closely related to the following system of PPDEs:
\bea
\label{PPDE}
\left\{\ba{lll}
\dis \pa_t u_i + {1\over 2} \si\si^\top(t, \o, u, \pa_\o u) : \pa^2_{\o\o} u_i  + f_i(t, \o, u, \pa_\o u) =0; \\
\dis u(T, \o) = g(\o),
\ea\right.\q i=1,\cds, d_2.
\eea

\subsection{Nonlinear Feynman-Kac formula}

The following result is an extension of the four step scheme of Ma, Protter, \& Yong \cite{MPY}.
\begin{thm}
\label{thm-FK}
Let Assumption \ref{assum-basic} hold, and $b, \si$ be  uniformly Lipschitz continuous in $(\bx, y,z)$. Assume  PPDE \reff{PPDE} has a classical solution $u\in C^{1,2}(\L; \dbR^{d_2})$ such that  $\pa_\o u$, $\pa^2_{\o\o} u$ are bounded and $u, \pa_\o u$ are uniformly Lipschitz continuous in $\o$. Then FBSDE  \reff{FBSDE} admits a strong solution and it holds that
\bea
\label{YZ}
 Y_t = u(t,X_\cd),\q Z_t =\pa_\o u(t,X_\cd).
\eea
Moreover, the solution is unique (in law) among all weak solutions. 
\end{thm}

\proof {\it Existence.} Set
\beaa
 \dis \tilde{b}(t,\o):=b(t,\o,u(t,\o),\pa_\o u(t,\o)),\q
  \tilde{\sigma}(t,\o):=\sigma(t,\o,u(t,\o),\pa_\o u(t,\o)).
 \eeaa
Under our conditions, both $\tilde{b}(t,\o)$ and $\tilde{\sigma}(t,\o)$ are bounded and are uniformly Lipschitz continuous in $\o$.
Thus, for any $x\in \dbR^{d_1}$, the following forward SDE
\bea
\label{SDE-X}
X_t= x+ \int_0^t \tilde{b}(s,X_\cd)ds+\int_0^t \tilde{\sigma}(s,X_\cd) dB_s,\q t\in [0,T],
\eea
has a (unique) strong solution.  Define $(Y, Z)$ by \reff{YZ} and $N_t:= 0$.  By applying functional It\^{o}'s formula \reff{FIto}, we can easily verify \reff{FBSDE}, hence $(X,Y, Z)$ is a strong solution of \reff{FBSDE}.

{\it Uniqueness.}  For notational simplicity let's assume $d_2=1$. The multidimensional case can be proved similarly without any significant difficulty.  Let  $(B, \Th, N, \dbP)$  be an arbitrary weak solution of \reff{FBSDE}.  We first claim that \reff{YZ} holds. Indeed, denote
\beaa
  \tilde{Y}_t=u(t,X_\cd),\q \tilde{Z}_t=\pa_\o u(t,X_\cd),\q \D Y_t := \tilde Y_t - Y_t,\q \D Z_t := \tilde Z_t - Z_t. 
\eeaa
Applying functional It\^{o} formula \reff{FIto} on $u(t,X_\cd)$ and recalling  \reff{PPDE}, we have:
\beaa
&&d \D Y_t  =  du(t,X_\cd) + f(t,X_\cd, Y_t, Z_t)dt - Z_t  dX_t + dN_t\\
&&=  \Big[\pa_t u(t,X_\cd)+{1\over 2}\pa^2_{\o\o}u(t,X_\cd):\sigma\si^\top(t, X_\cd, Y_t, Z_t) + f(t,X_\cd, Y_t, Z_t)\Big]dt +\D Z_t dX_t + dN_t\\
 &&=  -\Big[{1\over 2}\pa^2_{\o\o}u(t,X_\cd): \sigma\si^\top(t,X_\cd, \tilde Y_t, \tilde Z_t) + f(t,X_\cd, \tilde Y_t, \tilde Z_t)\Big] dt\\
 &&\qq +\Big[{1\over 2}\pa^2_{\o\o}u(t,X_\cd): \sigma\si^\top(t, X_\cd, Y_t, Z_t) + f(t,X_\cd, Y_t, Z_t)\Big]dt +\D Z_t dX_t +dN_t\\
&&=[\a_t \D Y_t + \b_t \D Z_t]dt + \D Z_t \si(t, X_\cd, Y_t, Z_t) d B_t + dN_t,
\eeaa
where $\a, \b$ are bounded. Note that $\D Y_T=0$. Applying It\^{o} formula on $|\D Y_t|^2$ and recalling Assumption \ref{assum-basic} (iii) we have
\beaa
&&\dbE\Big[|\D Y_t|^2 + c^2_0 \int_t^T |\D Z_s|^2 ds + tr(\la N\ra_T-\la N\ra_t)\Big] \\
&\le& \dbE\Big[|\D Y_t|^2 +  \Big|\int_t^T \D Z_s \si(s, X_\cd, Y_s, Z_s) d B_s\Big|^2 + \tr(\la N\ra_T-\la N\ra_t)\Big] \\
&=& \dbE\Big[\int_t^T 2\D Y_s  [\a_s \D Y_s + \b_s \D Z_s]ds\Big] \le  \dbE\Big[\int_t^T [C|\D Y_s|^2 + {c^2_0\over 2}|\D Z_s|^2]ds\Big] 
\eeaa
Then by the standard BSDE arguments we have $|\D Y| =|\D Z|=0$. This proves \reff{YZ}.

Now plug \reff{YZ} into the forward SDE of \reff{FBSDE}, we see that $X$ has to satisfy the SDE \reff{SDE-X}. By the uniqueness of \reff{SDE-X} we see that $X$ is unique, which, together with \reff{YZ}, implies further the uniqueness of $\Th$, hence that of $N$.
\qed

\section{Wellposedness for Markovian weak FBSDEs}
\label{sect-viscosity}
\setcounter{equation}{0}

We now turn to weak solutions. We shall follow the approach in Ma, Zhang, \& Zheng \cite{MZZ} and Ma \& Zhang \cite{MZ}. Our approach will rely heavily on viscosity solutions as well as the a priori estimates for the related PDE. We remark that all the results can be easily extended to path dependent case provided that  the corresponding estimates can be established for PPDEs, which however are not available in the literature and are in general challenging.  We thus restrict to Markovian case, and for the purpose of viscosity theory, we assume $d_2 = 1$.  Moreover, by Proposition \ref{prop-b}, we may assume without loss of generality that $b=0$. That is, our  objective of this section is the following  weak FBSDE:
 \bea
 \label{FBSDE1}
 \left\{\ba{lll}
 \dis X_t = x + \int_0^t \si(s, \Th_s) dB_s;\\
 \dis Y_t = g(X_T) + \int_t^T f(s, \Th_s) ds - \int_t^T Z_s dX_s + N_T-N_t,
 \ea\right. \dbP_0\mbox{-a.s.}
 \eea
  In this case the PPDE \reff{PPDE} becomes a standard  quasi-linear PDE: 
  \bea
\label{PDE}
\cL u(t,x) := \pa_t u(t,x) + {1\over 2} \si\si^\top (t,x,u,\pa_x u)  : \pa_{xx}^2 u  + f(t,x,u, \pa_x u)=0,~ u(T,x) = g(x),
\eea
extending  \reff{PDE-weak} to multidimensional case, and \reff{YZ} becomes
\bea
\label{YZ1}
Y_t = u(t, X_t), \q Z_t = \pa_x u(t, X_t),\q\dbP\mbox{-a.s.}
\eea

By Proposition \ref{prop-equiv}, throughout this section, we shall assume 
 \begin{assum}
 \label{assam-basic} (i) $d:= d_0=d_1$, $d_2=1$,  and $\si, f, g$ are state dependent;
 
 (ii) $\si, f(t,x,0,0), g$  are bounded by $C_0$, and $\si, f$ are continuous in $t$;
 
 (iii) $\si, f,  g$ are uniformly Lipschitz continuous in $(x,y,z)$ with Lipschitz constant $L$;
 
 (iv) $\si$ is symmetric and is uniformly nondegenerate: $\si  \ge c_0 I_d$ for some $c_0>0$;
 
 (v) Either $|\si(t,x,y,z_1) - \si(t,x,y, z_2)| \le {C_0\over 1+|z_1|} |z_1-z_2|$, or $d=1$. 
 \end{assum}
\no We emphasize again that, by Propositions \ref{prop-equiv} and \ref{prop-b}, we may allow $d_0 \neq d_1$ and \reff{FBSDE1} may depend on $b(t, X_\cd, Y_t, Z_t)$ as well.  Throughout this section, we use a generic constant $C>0$ which depends only on $T$ and  $C_0, c_0, L, d$ in Assumption \ref{assam-basic}.
 
 Under the above assumption, we have the following regularity results for the PDE \reff{PDE}. The arguments are mainly from Ladyzenskaja, Solonnikov \& Uralceva \cite{LSU}, and we sketch a proof in Appendix.
 \begin{thm}
 \label{thm-PDE}
 Let Assumption \ref{assam-basic} hold. Assume further that $\si, f, g$ are smooth with bounded derivatives. Then 
 
 (i) PDE \reff{PDE} has a classical solution $u\in C_b^{1, 2}([0,T]\times \dbR^d)$. 
 
(ii) There exists a constant $ \a>0$, depending only on $T$ and $C_0, c_0, L, d$ in Assumption \ref{assam-basic} , but not on the derivatives of $\si, f, g$, such that, for any $\d>0$,
 \bea
 \label{PDEest}
 \left.\ba{c}
 |u|\le C, \q |\pa_x u|\le C,\q  |u(t_1, x)- u(t_2,x)| \le C|t_1-t_2|^{1\over 2};\\
   |\pa_x u(t_1, x_1) -\pa_x u(t_2, x_2)|  \le C_\d\big[|x_1-x_2|^\a+|t_1-t_2|^{\a\over 2}\big] \q 0\le t_1 < t_2 \le T-\d.
 \ea\right.
 \eea
 where $C_\d$ may depend on  $\d$ as well.
 
 (iii) There exists a constant $C_g$, which depends on the same parameters $T, C_0, c_0, L, d$, as well as $\|\pa_{xx} g\|_\infty$, such that $|\pa_{xx} u|\le C_g$.
 \end{thm}

\subsection{Existence}

\begin{thm}
\label{thm-existence}
Let Assumption \ref{assam-basic} hold. Then FBSDE  \reff{FBSDE} admits a bounded semi-strong solution $\Th$, and \reff{YZ1} holds where $u$ is a viscosity solution of PDE \reff{PDE}.
\end{thm}
\proof Let $(\si_n, f_n, g_n)$ be a smooth mollifier of $(\si, f, g)$ such that they satisfy Assumption \ref{assam-basic}  uniformly. Applying Theorem \ref{thm-PDE}, let $u_n$ be the classical solution to PDE \reff{PDE} with coefficients $(\si_n, f_n, g_n)$, and then $\{u_n\}_{n\ge 1}$ satisfy \reff{PDEest} uniformly, uniformly in $n$. Applying the Arzela-Ascoli theorem, possibly along a subsequence,  $u_n$ converges to a function $u$ such that $u$ satisfies \reff{PDEest} and the convergence of $(u_n, \pa_x u_n)$ to $(u, \pa_x u)$ is uniform.  In particular, by the stability of viscosity solutions we see that $u$ is a viscosity solution of PDE \reff{PDE}.

Next, by Proposition \ref{prop-equiv} and Theorem \ref{thm-FK} the martingale problem \reff{FBSDE1}  with coefficients $(\si_n, f_n, g_n)$ has a solution $(\dbP_n, Z^n)$ such that $Y_t = u_n(t, X_t), Z^n_t = \pa_x u_n(t, X_t)$, $\dbP_n$-a.s. By Zheng \cite{Zheng}, possibly along a subsequence, we see that $\dbP_n$ converges to some $\dbP$ weakly. By the uniform convergence of $(u_n, \pa_x u_n)$, we have $Z^n_t \to Z_t$ uniformly, and  \reff{YZ1} holds. Moreover, it follows from \reff{PDEest} that $(Y, Z)$ are bounded. Finally, by the uniform convergence, it is straightforward  to verify that $(\dbP, Z)$ solves the martingale problem \reff{FBSDE1}  with coefficients $(\si, f,  g)$.   
\qed

\subsection{Nodal sets}

For any $t\in [0, T]$, we first extend Definition \ref{defn-weak} to  interval $[t, T]$.

\begin{defn}
\label{defn-weak2}
Let $(t,x,y)\in [0, T]\times \dbR^d\times \dbR$. We say $(B, \Th, N, \dbP)$ is a weak solution of FBSDE \reff{FBSDE1} at $(t,x,y)$ if they are processes on $[t, T]$ satisfying the requirements in Definition \ref{defn-weak} on $[t, T]$ and  $\dbP(X_t = x) = \dbP(Y_t = y) =1$. Define semi-strong solution, strong solution, and martingale problem at $(t,x,y)$ in an obvious sense.
\end{defn}
  We next define the nodal sets.

\begin{defn}
\label{defn-nodal}
(i) For  $(t,x, y)\in [0, T]\times \dbR^d\times \dbR$, let $\cO(t,x,y)$ denote the space of weak solutions   of \reff{FBSDE1} at $(t,x,y)$.

(ii) $O := \{(t,x,y):  \cO(t,x,y)\neq \emptyset\}$.
\end{defn}

By Theorem \ref{thm-existence}, $(t,x, u(t,x))\in O$ for any $(t,x) \in [0, T]\times \dbR^d$,  where $u$ is the viscosity solution of PDE \reff{PDE} in Theorem \ref{thm-existence}.   We remark that a priori the measurability of $O$ is not clear. Nevertheless,  let $\ol O$ denote the closure of $O$, and define
 \bea
\label{baru}
\ul u(t, x) := \inf\{y: (t,x,y)\in \ol O\},\q \ol u(t,x) := \sup\{y: (t,x,y)\in \ol O\}.
\eea
Then $\ul u$ and $\ol u$ are Borel measurable and $\ul u \le u\le \ol u$.
 
\begin{prop}
\label{prop-baru}
Let Assumption \ref{assam-basic} hold. Then

(i)  $\ol u$ and $\ul u$ are bounded;

(ii) $\ol u$ is  upper semi-continuous and $\ul u$ is lower semi-continuous;

(iii) $\ol u(T,x) = \ul u(T,x) = g(x)$.
\end{prop}
\proof (i) For any $(t,x, y)\in O$ with corresponding weak solution $(\Th, N, \dbP)$, we have
\beaa
Y_s = g(X_T) + \int_s^T f(r, \Th_r) dr - \int_t^T Z_r dX_r + N_T-N_t.
\eeaa
Since $g$ and $f(t, x,,0,0)$ are bounded by $C_0$ and $f$ is uniformly Lipschitz continuous in $(y,z)$, it follows from standard BSDE arguments that 
\bea
\label{Ybound}
\dbE^{\dbP}\Big[\sup_{t\le s\le T}[|Y_s|^2+|N_s|^2]+ \int_t^T |Z_s|^2 ds\Big] \le C.
\eea
  In particular, $|y| = |Y_t|\le C$. This implies $|\ul u|, |\ol u|\le C$.

Since $\ol O$ is closed, (ii) is a direct consequence of the definitions of $\ol u$, $\ul u$. To see (iii), let $(T,x,y)\in \ol O$. By definition there exist  $(t_n, x_n, y_n)\in O$ such that $t_n \uparrow T$ and $(x_n, y_n) \to (x, y)$. Let $(B^n, \Th^n, \dbP^n)$ be a weak solution at $(t_n, x_n, y_n)$. Then
\beaa
|y_n - g(x_n)|^2 &=& \Big|\dbE^{\dbP_n} \Big[g(X^n_T) + \int_{t_n}^T f(s, \Th^n_s)  ds\Big] - g(x_n)\Big|^2\\
&\le& C \dbE^{\dbP_n}\Big[|X^n_T - x_n|^2  + (T-t_n) \int_{t_n}^T [1+ |Y_s|^2 + |Z^n_s|^2]ds\Big]  \\
&\le& C \dbE^{\dbP_n}\Big[ \int_{t_n}^T |\si(s, \Th^n_s) |^2 ds\Big]  +C (T-t_n) \le C (T-t_n),
\eeaa
thanks to \reff{Ybound}. Send $n\to \infty$, we see that $y = g(x)$. This proves (iii).
\qed

We have the following result improving Theorem \ref{thm-existence}, which is not used in this paper but is nevertheless interesting in its own right.

\begin{thm}
\label{thm-existence2}
Let Assumption \ref{assam-basic} hold. Then $(t,x,y) \in O$ if and only if $y \in [\ul u(t,x), \ol u(t,x)]$. Moreover, for any $(t,x,y)\in O$, there exists a semi-strong solution $(B, \Th, N, \dbP)$ at $(t,x,y)$ such that $|Z|\le C$. 
\end{thm}
\proof It is clear that $(t,x,y) \in O$ implies $y \in [\ul u(t,x), \ol u(t,x)]$. Then it suffices to show that, for any  $y \in [\ul u(t,x), \ol u(t,x)]$, there exists a weak solution at $(t,x,y)$ such that $Z$ is bounded. We proceed in two steps.

{\it Step 1.} For any $n\ge 1$, let $\si_n, f_n, g_n$ be smooth mollifiers of $\si, f, g$ such that
\bea
\label{mollifier}
|\si_n - \si|\le \e_n,\q |f_n - f|\le {1\over n},\q |g_n - g|\le {1\over n},
\eea
for some small $\e_n>0$ which will be specified later. Denote 
\beaa
\ol f_n := f_n + {2\over n},\q \ul f_n := f_n - {2\over n},\q \ol g_n := g_n + {1\over n},\q \ul g_n := g_n - {1\over n}.
\eeaa
By Theorem \ref{thm-PDE}, the PDE \reff{PDE} with coefficients $(\si_n, \ol f_n, \ol g_n)$ (resp. $(\si_n, \ul f_n, \ul g_n)$) has a classical solution $\ol u_n$ (resp. $\ul u_n$). We claim that, for any $(t,x,y)\in O$ and any $n$,
\bea
\label{comparison}
\ul u_n(t,x) \le y \le \ol u_n(t,x).
\eea

Without loss of generality we will prove only the right inequality at $t=0$.  We shall follow similar arguments as in Theorem  \ref{thm-FK}. Let $(B, \Th, N, \dbP)$ be a weak solution to FBSDE \reff{FBSDE1}   at $(0,x,y)$ with coefficients $(\si, f, g)$. Fix $n$ and denote
\beaa
\tilde Y_t := \ol u_n(t, X_t),~ \tilde Z_t := \pa_x \ol u_n(t, X_t),~ \tilde \Th := (X, \tilde Y, \tilde Z),~  \D Y_t := \tilde Y_t - Y_t,~ \D Z_t := \tilde Z_t - Z_t.
\eeaa
Apply It\^{o} formula, we have
\beaa
d\D Y_t &=& \Big[\pa_t \ol u_n + {1\over 2} \pa^2_{xx} \ol u_n : \si^2(t, \Th_t) + f(t, \Th_t)\Big]dt + \D Z_t dX_t + dN_t\\
&=& \Big[ {1\over 2} \pa^2_{xx} \ol u_n : [\si^2(t, \Th_t) -\si^2_n(t,  \tilde \Th_t)] + [f(t, \Th_t)  -\ol f_n(t, \tilde \Th_t)]\Big]dt  + \D Z_t dX_t + dN_t.
\eeaa
By Theorem \ref{thm-PDE} (iii), there exists a constant $C_n$, which is independent of $\e_n$, such that $|\pa^2_{xx} \ol u_n|\le C_n$. Note that $\ol f_n - f = f_n +{2\over n}  - f \ge {1\over n}$ and $|\si|\le C_0$.  Then, for $\e_n\le {1\over n C_0C_n}$, we have
\beaa
d\D Y_t &\le& \Big[ {1\over 2} \pa^2_{xx} \ol u_n : [\si^2(t, \Th_t) -  \si^2(t,  \tilde \Th_t)] + [f(t, \Th_t) -f(t, \tilde \Th_t)]\Big]dt  + \D Z_t dX_t + dN_t\\
&=&  \Big[ \a^n_t \D Y_t + \b^n_t \D Z_t\Big]dt  + \D Z_t dX_t + dN_t,
\eeaa
where $|\a^n|, |\b^n|\le C_n$. Note further that $\D Y_T = \ol g_n(X_T) - g(X_T) = g_n(X_T) +{1\over n} - g(X_T)\ge 0$. It is clear that $\D Y_0 \ge 0$. This implies  $0\le \tilde Y_0 - Y_0 = \ol u_n(0, x) - y$, proving \reff{comparison}.

 {\it Step 2.} Let $y \in [\ul u(t,x), \ol u(t,x)]$. There exist $(\ul t_m, \ul x_m, \ul y_m)\in O$ and  $(\ol t_m, \ol x_m, \ol y_m)\in O$ such that $(\ul t_m, \ul x_m, \ul y_m)\to (t,x, \ul u(t,x))$ and  $(\ol t_m, \ol x_m, \ol y_m)\to (t,x,\ol u(t,x))$.  Then, by \reff{comparison},
 \beaa
 \ul u_n(\ul t_m, \ul x_m)  \le \ul y_m,\q \ol y_m \le \ol u_n(\ol t_m, \ol x_m),\q \mbox{for all}~m, n. 
\eeaa
Send $m\to \infty$, we obtain
\bea
\label{yolu}
 \ul u_n(t, x)  \le \ul u(t,x) \le y\le  \ol u(t,x) \le \ol u_n(t, x),\q \mbox{for all}~ n. 
 \eea
 For any $n\ge 1$ and $\a\in [0, 1]$, denote $\f^\a_n := \a \ol \f_n + [1-\a] \ul \f_n$ for $\f = f, g$, and let $u^\a_n$ be the classical solution of PDE \reff{PDE} with coefficients $(\si_n, f^\a_n, g^\a_n)$.  By the arguments in Theorem \ref{thm-PDE}, it is clear that the mapping $\a\mapsto u^\a_n(0,x)$ is continuous. Since $u^0_n(t,x) = \ul u_n(t,x) \le y \le \ol u_n(t,x) = u^1_n(t,x)$, there exists $\a_n \in [0,1]$ such that $u^{\a_n}_n(t,x) = y$. For each $n\ge 1$,  by Proposition \ref{prop-equiv} and Theorem \ref{thm-FK}   the martingale problem \reff{FBSDE1} at $(t,x, y)$ with coefficients $(\si_n, f^{\a_n}_n, g^{\a_n}_n)$ has a solution $(\dbP^n, Z^n)$ such that $Y_s = u_n^{\a_n}(s, X_s), Z^n_s = \pa_x u_n^{\a_n}(s, X_s)$, $t\le s\le T$, $\dbP^n$-a.s.  Now following the arguments in Theorem \ref{thm-existence} we  see that, possibly following a subsequence,  $\dbP^n \to \dbP$, $Z^n\to Z$, $u^{\a_n}_n \to u$, where $(\dbP, Z)$ is a solution to  the martingale problem \reff{FBSDE1} at $(t,x,y)$ with coefficients $(\si, f, g)$ and $u$ is a viscosity solution to PDE \reff{PDE} with coefficients $(\si, f, g)$.  It is clear that $|Z_s| = |\pa_x u(s, X_s)|\le C$, $\dbP$-a.s.
 \qed

\subsection{Uniqueness}
\begin{thm}
\label{thm-viscosity}
Let Assumption \ref{assam-basic} hold.  Then $\ol u$ (resp. $\ul u$) is a viscosity subsolution (resp. supersolution)  of PDE \reff{PDE}. 
\end{thm}
\proof We shall prove the result only for $\ol u$. The result for $\ul u$ can be proved similarly.

Fix $(t_0, x_0)\in [0, T)\times \dbR^d$ and denote $y_0 := \ol u(t_0, x_0)$. Let $\f \in C^{1,2}_b([0, T]\times \dbR^d)$ be a test function at $(t_0, x_0)$, namely
\beaa
[\f- \ol u](t_0, x_0) =  0 = \inf_{(t,x) \in [0, T]\times \dbR^d} [\f-\ol u](t,x).
\eeaa
Let $(t_n, x_n, y_n)\in O$ such that $(t_n, x_n, y_n) \to (t_0, x_0, y_0)$, and  $(\dbP^n, Z^n)$ a weak solution to the martingale problem \reff{FBSDE1} at $(t_n, x_n, y_n)$. Define $N^n$ as in \reff{BN}. By using regular conditional probability distribution, it is clear that  $(t, X_t, Y_t) \in O$, $\dbP^n$-a.s. for $t_n \le t\le T$. Then, by the definition of $\ol u$, we have $Y_t \le \ol u(t, X_t) \le \f(t, X_t)$. 

Now denote 
\beaa
&\D Y_t := \f(t, X_t) - Y_t \ge 0,\q \D Z^n_t := \pa_x \f(t, X_t) - Z^n_t,&\\
& \Th^n_t := (X_t, Y_t, Z^n_t),\q \tilde \Th_t := (X_t, \f(t, X_t), \pa_x \f(t, X_t)).&
\eeaa
Applying It\^{o} formula we have, under $\dbP^n$,
\beaa
d\D Y_t &=&  \Big[\pa_t \f(t, X_t) + {1\over 2} \pa^2_{xx} \f(t, X_t) : \si^2(t, \Th^n_t)  dt + f(t, \Th^n_t)\Big] dt + \D Z^n_t dX_t + dN^n_t\\
&=& \Big[\cL \f(t, X_t) +  {1\over 2} \pa^2_{xx} \f(t, X_t) : [\si^2(t, \Th^n_t) - \si^2(t, \tilde \Th_t)] + [f(t, \Th^n_t) - f(t, \tilde \Th_t)] \Big] dt \\
&&+ \D Z^n_tdX_t + dN^n_t\\
&=&  \Big[\cL \f(t, X_t) -  \a^n_t\D Y_t -  \D Z^n_t \si(t, \Th^n_t) \b^n_t\Big] dt + \D Z^n_tdX_t + dN^n_t,
\eeaa
where $|\a^n|, |\b^n|\le C$. Denote
\beaa
\G^n_t := \exp\Big( \int_{t_n}^t \b^n_s\cd \si^{-1}(s, \Th^n_s) dX_s + \int_{t_n}^t [\a^n_s - {1\over 2} |\b^n_s|^2] ds\Big).
\eeaa
Then
\beaa
d[\G^n_t \D Y_t] =  \G^n_t \cL\f(t, X_t) dt +  \G^n_t[\D Z^n_t + \D Y_t [\b^n_s]^\top \si^{-1}(s, \Th^n_s)] dX_t+ \G^n_t dN^n_t.
\eeaa
Thus, for any $\d>0$ small, 
\beaa
0 &\le& \dbE^{\dbP_n}[\G^n_{t+\d}\D Y_{t_n +\d}] =  \dbE^{\dbP_n}\Big[  \G^n_{t_n} \D Y_{t_n} +\int_{t_n}^{t_n+\d} \G^n_t \cL\f(t, X_t) dt\Big]\\
&=& \f(t_n, x_n)-y_n +  \cL \f(t_n, x_n) \d  +   \dbE^{\dbP_n}\Big[\int_{t_n}^{t_n+\d}[ \G^n_t \cL\f(t, X_t)  - \G^n_{t_n} \cL\f(t_n, X_{t_n}]dt\Big].
\eeaa
Note that $\cL\f$ is uniformly continuous, and since $\si$ is bounded, one can easily show that
\beaa
  \dbE^{\dbP_n}\Big[\big| \G^n_t \cL\f(t, X_t)  - \G^n_{t_n} \cL\f(t_n, X_{t_n}\big|\Big] \le \rho(\d),\q t_n \le t\le t_n+\d,
 \eeaa
 for some modulus of continuity function $\rho$. Then
 \beaa
0\le    \f(t_n, x_n)-y_n + \cL \f(t_n, x_n) \d   + \d \rho(\d).
  \eeaa
  Send $n\to \infty$, we have
  \beaa
 0\le  \cL \f(t_0, x_0) \d + \d \rho(\d).
  \eeaa
Divide both sides by $\d$ and then send $\d \to 0$, we obtain $\cL\f(t_0, x_0) \ge 0$.
\qed  

We remark that, in the case that $\si$ is independent of $z$, \cite{MZ} and \cite{MZZ} established similar results without requiring the uniform Lipschitz continuity of the coefficients, and thus the arguments there are more involved.  

Our final result relies on the comparison principle for viscosity solutions of PDEs, for which we refer to the classical reference Crandall, Ishii, \& Lions \cite{CIL}. We say a PDE satisfies the comparison principle for viscosity solutions if: {\it for any upper semi-continuous viscosity subsolution  $u_1$ and any lower semi-continuous viscosity supersolution $u_2$ with $u_1(T,\cd) \le u_2(T, \cd)$, we have $u_1\le u_2$.}
\begin{thm}
\label{thm-uniqueness} 
Let Assumption \ref{assam-basic} hold. Assume further that the comparison principle for the viscosity solutions of PDE \reff{PDE} holds true. Then the weak solution to FBSDE \reff{FBSDE1} is unique (in law).
\end{thm}
\proof First, by the comparison principle, it follows from Theorem \ref{thm-viscosity} that $\ol u = \ul u = u$, where $u$ is the unique viscosity solution of the PDE \reff{PDE} satisfying \reff{PDEest}. Now let $(B, \Th, N, \dbP)$ be an arbitrary weak solution of FBSDE \reff{FBSDE1}. Since $(t, X_t, Y_t)\in O$, $\dbP$-a.s., then $Y_t = u(t, X_t)$, $\dbP$-a.s.

Next, for any $\d>0$, $0< t\le T-\d$, and any partition $0=t_0<\cds<t_n = t$ with $ t_{i+1}-t_i = h := {t\over n}$, by \reff{PDEest} we have
\beaa
&&\Big|\sum_{i=0}^{n-1} [Y_{t_{i+1}} - Y_{t_i}][X_{t_{i+1}}-X_{t_i}] - \int_0^t \pa_x u(s, X_s) d \la X\ra_s\Big|\\
&=&\Big|\sum_{i=0}^{n-1} [u(t_{i+1}, X_{t_{i+1}}) - u(t_i, X_{t_i})][X_{t_{i+1}}-X_{t_i}] - \int_0^t \pa_x u(s, X_s) d \la X\ra_s\Big|\\
&\le&\Big|\sum_{i=0}^{n-1}  [u(t_{i+1}, X_{t_{i}}) - u(t_{i}, X_{t_i})][X_{t_{i+1}}-X_{t_i}] \Big| \\
&&+ \Big|\sum_{i=0}^{n-1}\Big[ [u(t_{i+1}, X_{t_{i+1}}) - u(t_{i+1}, X_{t_i})][X_{t_{i+1}}-X_{t_i}] -  \pa_x u(t_{i+1}, X_{t_i})  [\la X\ra_{t_{i+1}}-\la X\ra_{t_i}]\Big] \Big|\\
&&+\sum_{i=0}^{n-1} \Big|\int_{t_i}^{t_{i+1}} [\pa_x u(s, X_s) - \pa_x u(t_{i+1}, X_{t_i})]d \la X\ra_s\Big|\\
&\le& \Big|\sum_{i=0}^{n-1}  [u(t_{i+1}, X_{t_{i}}) - u(t_{i}, X_{t_i})][X_{t_{i+1}}-X_{t_i}] \Big| \\
&&+ \Big|\sum_{i=0}^{n-1} \pa_x u(t_{i+1}, X_{t_i})\big[[X_{t_{i+1}}-X_{t_i}] ^\top [X_{t_{i+1}}-X_{t_i}] -[\la X\ra_{t_{i+1}}-\la X\ra_{t_i}] \big]\Big| \\
&&+C_\d\sum_{i=0}^{n-1}  |X_{t_{i+1}}-X_{t_i}|^{2+\a} +C_\d \sum_{i=0}^{n-1}\int_{t_i}^{t_{i+1}} [h^{\a\over 2} + |X_{t_{i+1}} - X_s|^\a] ds.
\eeaa
Since $X_t-X_s = \int_s^t \si(r, \Th_r) dB_r$ and $\si$ is bounded, one can easily show that
\beaa
\dbE^\dbP[|X_t - X_s|^p] \le C_p |t-s|^{p\over 2}.
\eeaa
Moreover, by the martingale property of $X$, we have
\beaa
&& \dbE^\dbP\Big[ \Big|\sum_{i=0}^{n-1}  [u(t_{i+1}, X_{t_{i}}) - u(t_{i}, X_{t_i})][X_{t_{i+1}}-X_{t_i}] \Big|\Big]^2\\
&=& \dbE^\dbP\Big[\sum_{i=0}^{n-1} \Big| [u(t_{i+1}, X_{t_{i}}) - u(t_{i}, X_{t_i})][X_{t_{i+1}}-X_{t_i}] \Big|^2\Big]\\
&\le& C h \dbE^\dbP\Big[\sum_{i=0}^{n-1} |X_{t_{i+1}}-X_{t_i}|^2\Big] \le C h \sum_{i=0}^{n-1} h = C h;
\eeaa
and, applying It\^{o} formula, 
\beaa
&&\dbE^\dbP\Big[\Big|\sum_{i=0}^{n-1} \pa_x u(t_{i+1}, X_{t_i})\big[[X_{t_{i+1}}-X_{t_i}] ^\top [X_{t_{i+1}}-X_{t_i}] -[\la X\ra_{t_{i+1}}-\la X\ra_{t_i}] \big]\Big|^2\Big]\\
&=&\dbE^\dbP\Big[\Big|\sum_{i=0}^{n-1} \pa_x u(t_{i+1}, X_{t_i}) \int_{t_i}^{t_{i+1}} [X_s-X_{t_i}]^\top d X_s \Big|^2\Big]\\
&=&\dbE^\dbP\Big[\sum_{i=0}^{n-1} \Big|\pa_x u(t_{i+1}, X_{t_i}) \int_{t_i}^{t_{i+1}} [X_s-X_{t_i}]^\top d X_s \Big|^2\Big]\\
&\le& C\dbE^\dbP\Big[\sum_{i=0}^{n-1} \int_{t_i}^{t_{i+1}} |X_s-X_{t_i}|^2 d s \Big]\le \sum_{i=0}^{n-1} \int_{t_i}^{t_{i+1}} [s-t_i] d s \le Ch.
\eeaa
Then we have
\beaa
\dbE^\dbP\Big[\Big|\sum_{i=0}^{n-1} [Y_{t_{i+1}} - Y_{t_i}][X_{t_{i+1}}-X_{t_i}] - \int_0^t \pa_x u(s, X_s) d \la X\ra_s\Big|^2\Big]\le  Ch +C_\d h^\a   \le C_\d h^\a.
\eeaa
Send $n\to \infty$ and thus $h\to 0$, note that 
\beaa
\sum_{i=0}^{n-1} [Y_{t_{i+1}} - Y_{t_i}][X_{t_{i+1}}-X_{t_i}] \to \la Y, X\ra_t = \int_0^t Z_s d\la X\ra_s,\q \mbox{in}~ \dbL^2(\dbP),
\eeaa
then we have
\beaa
\int_0^t Z_s d\la X\ra_s = \int_0^t \pa_x u(s, X_s) d\la X\ra_s,\q \dbP\mbox{-a.s.},\q 0\le t\le T-\d.
\eeaa
Since $\si$ is nondegenerate and $t$ and $\d$ are arbitrary, we obtain
\beaa
Z_t = \pa_x u(t, X_t),\q dt \times d\dbP\mbox{-a.s. on}~[0, T) \times \O.
\eeaa
That is, \reff{YZ1} holds.

Now similar to the existence part of Theorem \ref{thm-FK}, denote
\beaa
\tilde \si(t,x) := \si(t, x, u(t,x),\pa_x u(t,x)).
\eeaa
Then $\tilde \si$ is H\"{o}lder continuous and $(B, X, \dbP)$ satisfies the SDE:
\beaa
X_t = x + \int_0^t \tilde \si(s, X_s) dB_s,\q \dbP\mbox{-a.s.}
\eeaa
By Stroock \& Varadahn \cite{SV}, the above SDE has a unique (in law) weak solution. This, together with \reff{YZ1}, implies the uniqueness (in law) of $(B, \Th, \dbP)$. Finally, by \reff{BN}, the joint law with $N$ is also unique.
\qed

\begin{rem}
\label{rem-target}
{\rm An alternative approach to prove the uniqueness is to consider the stochastic target problem, as in Soner and Touzi \cite{ST}. That is, in the spirit of \reff{baru0}, define
\beaa
&\dis \ol u(t,x) := \inf\Big\{y: \exists Z ~\mbox{such that}~ Y^{t,x,y,Z}_T \ge g(X^{t,x,Z}_T),~\dbP_0\mbox{-a.s.}\Big\},\q\mbox{where}&\\
&\dis X^{t,x,y,Z}_s = x + \int_t^s \si(r, X^{t,x,y,Z}_r) dB_r,&\\
&\dis Y^{t,x,y,Z}_s = y - \int_t^s f(r, X^{t,x,y,Z}_r, Y^{t,x,y,Z}_r, Z_r) dr + \int_t^s Z_r d X^{t,x,y,Z}_r,&
\eeaa
and define $\ul u$ similarly.  The idea is to prove that $\ol u$ and $\ul u$ are viscosity solutions of the PDE. However, there are technical difficulties in establishing the regularity and the dynamic programming principle for these functions. We shall leave this possible approach to future research.
\qed}
\end{rem}

\section{Appendix}
\label{sect-Appendix}

\subsection{Some counterexamples}
In this subsection we provide two counterexamples related to the control problems in Section \ref{sect-control}. In particular, they will show that the stochastic control problems in weak formulation have optimal controls, while the corresponding problems in strong formulation do not have optimal control. In the first example, we also show that the associated weak FBSDE has a weak solution, but no strong solution.

 \subsubsection{The case with drift control}
 In this case we shall consider an example with path dependence. We note that all the heuristic analysis in Section \ref{sect-control} can be easily extended to the path dependent case. We first recall a result due to  Tsirel'son \cite{Tsirelson}. 

\begin{lem}
\label{lem-Tsirel'son}
Let $t_n>0$, $n\ge 1$,  be strictly decreasing with $t_0 = T$ and $t_n \downarrow 0$, and $\th(x):=x-[x]$ where $[x]$ is the largest integer in $(-\infty, x]$. Define the non-curtailing functional $K$:
\bea
\label{K}
K(t, \bx) := \th(\frac{\bx(t_n)-\bx(t_{n+1})}{t_n-t_{n+1}}), ~\mbox{for}~ t\in [t_n,t_{n-1}), ~ \bx\in C([0, T]).
\eea
Then the following path dependent SDE  has no strong solution:
\bea
\label{Tsirel'son}
X_t= \int_0^t K(s,X_\cd)ds+B_t.
\eea
\end{lem}
\no We remark that $K$ is bounded and thus SDE \reff{Tsirel'son} has a unique (in law) weak solution, following the standard Girsanov Theorem. We also note that the above $K$ is discontinuous. When $K$ is state dependent, namely $K = K(t, X_t)$, the SDE could have a strong solution even when $K$ is discontinuous, see Cherny \& Engelbert \cite{CE}  and   Halidias \& Kloeden \cite{HK}  for some positive results.
 
 Our example considers the following setting, with $f$ depending on the paths of $X$:
\bea
\label{drift-setting}
b(t,\a) := \a,\q \si :=1,\q f(t,\bx,\a) := -{1\over 2}|\a- K(t,\bx)|^2,\q g := 0.
\eea

 \begin{eg}
\label{eg-drift} Let $K$ be defined in \reff{K}, and $\cA := \dbL^2(\dbF^B, \dbP_0)$.

(i) The   optimization problem in weak formulation has  an optimal control $\a^*_t := K(t, X_\cd)$:
 \bea
 \label{drift-weak}
 &\dis \bar V_0 := \sup_{\a\in \cA}  \bar V^\a_0 :=  \sup_{\a\in \cA}  \dbE^{\dbP^\a}\Big[-{1\over 2} \int_0^T |\a_t - K(t, X_\cd)|^2 dt\Big],&\\
  &\dis  \mbox{where}\q X_t := B_t,\q d\dbP^\a := M^\a_T d\dbP_0 := \exp\Big(\int_0^T \a_t dB_t - {1\over 2} \int_0^T|\a_s|^2 dt\Big)d\dbP_0.&\nonumber
 \eea

(ii) The  optimization problem in strong formulation has no optimal control:
\bea
\label{drift-strong}
&\dis V_0 := \sup_{\a\in\cA} V^\a_0,\q\mbox{where}&\\
&\dis X^\a_t := \int_0^t \a_s ds + B_t,\q  V^\a_0 = \dbE^{\dbP_0}\Big[-{1\over 2} \int_0^T |\a_t - K(t, X^\a_\cd)|^2 dt \Big].&\nonumber
\eea

\end{eg}
 \proof (i) Since $\bar V^\a_0\le 0$, it is obvious that $\bar V_0 \le 0$. Moreover, it is clear that $\bar V^{\a^*} =0$ for $\a^*_t := K(t, X_\cd) = K(t, B_\cd)$, then $\bar V_0=0$ with optimal control $\a^*$.
 
 (ii) For each $n$, denote $t_i := {iT\over n}$, $i=0,\cds, n$, and $\a^n_t := \sum_{i=1}^{n-1} {T\over n} \int_{t_{i-1}}^{t_i} K(s, B_\cd)ds  \1_{[t_i, t_{i+1})}$. Recall \reff{drift-weak} and note that
 \beaa
   \bar V^\a_0 =   \dbE^{\dbP_0}\Big[M^\a_T \big[-{1\over 2} \int_0^T |\a_t - K(t,  B_\cd)|^2 dt\big]\Big].
   \eeaa
 It is clear that 
 \bea
 \label{barVan}
 \lim_{n\to\infty} \bar V^{\a^n}_0 = \bar V^{\a^*}_0 = 0.
 \eea
 Since $\a^n$ is piecewise constant, then $\dbF^{B} = \dbF^{B^{\a^n}}$, and thus there exists a piecewise constant process $\tilde \a^n$ such that $\a^n_t(B_\cd) = \tilde \a^n_t(B^{{\a^n}}_\cd)$. That is,
 \beaa
 B_t = \int_0^t \tilde \a^n_s(B^{{\a^n}}_\cd) ds + B^{\a^n}_t,\q \dbP^{\a^n}\mbox{-a.s.}
 \eeaa
 Therefore, the $\dbP_0$-distribution of $(B, X^{\tilde \a^n}, \tilde \a^n(B))$ coincides with the $\dbP^{\a^n}$-distribution of $(B^{\a^n}, B, \a^n(B))$.  This implies that $V^{\tilde \a^n}_0 = \bar V^{\a^n}_0$. Then by \reff{barVan} we see that $V_0 \ge \lim_{n\to\infty} V^{\tilde \a^n}_0 =0$. On the other hand, it is obvious that $V_0 \le 0$. Then $V_0=0$.

 Now if \reff{drift-strong} has an optimal control $\tilde \a^*$, then $V^{\tilde\a^*}_0 = 0$ and thus $\tilde \a^*_t = K(t, X^{\tilde \a^*}_\cd)$, $\dbP_0$-a.s.  Thus $X^*:= X^{\tilde\a^*}$ satisfies SDE \reff{Tsirel'son}. 
 Since by definition $\tilde \a^*$ is $\dbF^B$-progressively measurable, we see that $X^*$ is also $\dbF^B$-progressively measurable, and hence $X^*$ is a strong solution of SDE \reff{Tsirel'son}. This contradicts with Lemma \ref{lem-Tsirel'son}.
 \qed 

\begin{rem}
\label{rem-drift}
{\rm By extending the arguments to this case, one can (formally) show that the weak FBSDE \reff{FBSDEcontrol5} and the equivalent one \reff{XbarY} becomes
   \bea
 \label{drift-FBSDE}
 \left\{\ba{lll}
 \dis X_t = \int_0^t  [ Z_s + K(s, X_\cd)] ds +B_t;\\
 \dis Y_t = \int_t^T [ {1\over 2}|Z_s|^2 +  K(s, X_\cd) Z_s] ds  - \int_t^T Z_s dX_s.
  \ea\right.
 \eea
This FBSDE has a weak solution: $Y=Z=0$ and $X$ is the weak solution to SDE \reff{Tsirel'son}. However, it does not have a strong solution such that $\int_0^t Z_s dB_s$ is a BMO martingale. We refer to Zhang \cite{Zhang-book} Chapter 7 for BMO martingales. Indeed, if there is such a solution, then by \reff{drift-FBSDE} we immediately have
\beaa
Y_t =  -{1\over 2} \int_t^T |Z_s|^2  ds  - \int_t^T Z_s dB_s.
\eeaa
This implies that $Y=Z=0$. Then $X$ has to be a strong solution of SDE \reff{Tsirel'son}, contradicting with Lemma \ref{lem-Tsirel'son}.
\qed}
\end{rem}

  \subsubsection{The case with diffusion control}

We first recall a result due to Barlow \cite{Barlow}.  Recall the function $\th(x)$ in Lemma \ref{lem-Tsirel'son}. 
\begin{lem}
\label{lem-Barlow}
Let $\frac{\sqrt{2}}{2}<\l<1$ and define
\bea
\label{Barlow-sigma}
\si_0(x) := 1+  \sum^{\infty}_{n=0} \l^n \eta\big(\th(2^n x)\big),\q\mbox{where}\q \eta(x):= x \bold{1}_{[0,\frac{1}{2})}(x)+(1-x) \bold{1}_{[\frac{1}{2},1)}(x).
\eea
Then the following SDE has a unique weak solution but no strong solution:
\bea
\label{Barlow-SDE}
X_t=  \int_0^t \si_0(X_s) dB_s,\q\dbP_0\mbox{-a.s.} 
\eea
\end{lem}
\proof We first note that, although $\th$ is discontinuous at integers, $\eta \circ \th$ is actually Lipschitz continuous and periodic. Then $\si_0$ is uniformly continuous, and clearly $\si_0\ge 1$. Thus it follows from Stroock \& Varadahn \cite{SV} that \reff{Barlow-SDE} has a unique weak solution.

On the other hand, one may verify that $\si_0$ satisfies the hypotheses in \cite{Barlow}  Theorem 1.3 with $\alpha=\beta=-\ln \l\slash \ln 2$.  Then we see that  \reff{Barlow-SDE} has no strong solution.
\qed

The next example considers the following setting with diffusion control:
\bea
\label{diffusion-setting}
b := 0,~ \si(t,\a) :=\a,~ f(t,x,\a) := -{1\over 4}[|\a|^4 + |\si_0(x)|^4],~ g(x):=\int^x_0\int^\l_0 [\si_0(r)]^2 dr d\l.
\eea
Note that in this case we need the weak formulation for diffusion control problems. We refer to Zhang \cite{Zhang-book} Chapter 9 for details.
\begin{eg}
\label{eg-Barlow} Consider \reff{Barlow-sigma}  and \reff{diffusion-setting} with $\l={3 \over 4}$, and let the control set $A := [1,2]$.

(i) The optimization problem in weak formulation has optimal control $\a^*_t := \si_0(X_t) $:
 \bea
 \label{diffusion-weak}
 &\dis \bar V_0 := \sup_{\a\in \cA}  \bar V^\a_0 :=  \sup_{\a\in \cA}  \dbE^{\dbP^\a}\Big[g(X_T) + \int_0^T f(t, X_t, \a_t)  dt\Big],&\\
  &\dis  \mbox{where $\dbP^\a$ is a weak solution of SDE:}~ X_t = \int_0^t \a_s(X_\cd) dB_s. &\nonumber
 \eea

(i) The optimization problem in strong formulation has no optimal control:
 \bea
\label{diffusion-strong}
&\dis V_0 := \sup_{\a\in \cA} V^\a_0:= \sup_{\a\in \cA} \dbE^{\dbP_0}\Big[g(X^\a_T) + \int_0^T f(t, X^\a_t, \a_t) dt\Big].&\\
  &\dis  \mbox{where}\q X^\a_t := \int_0^t\a_s(B_\cd) dB_s.&\nonumber
  \eea

\end{eg}
 \proof 
  (i) By standard literature,  $\bar V_0 = u(0,0)$, where $u$ satisfies the HJB equation:
 \bea
 \label{diffusion-PDE}
 \pa_t u +\sup_{\a\in [1,2]} \Big[{1\over 2}\a^2\pa^2_{xx} u - {1\over 4}\a^4 \Big]- {1\over 4}|\si_0(x)|^4=0,\q u(T,x) = g(x).
 \eea
Note that the above PDE has a classical solution $u(t,x)=g(x)$. Then $\bar V_0=g(0) =0$.  On the other hand, let $\a^*_t (X_\cd) := \si_0(X_t)$ and $\dbP^*:= \dbP^{\a^*}$ be the (unique) weak solution of SDE \reff{Barlow-SDE}. Denote $Y_t := g(X_t) + \int_0^t f(s, X_s, \a^*_s)ds$ and note that $g''(x) = |\si_0(x)|^2$. Then applying It\^{o} formula we have
\beaa
d Y_t &=& \Big[{1\over 2} g''(X_t) |\si_0(X_t)|^2 + f(t, X_t, \a^*_t)\Big] dt + g'(X_t) \si_0(X_t) dB_t\\
&=&  \Big[{1\over 2} |\si_0(X_t)|^4- {1\over 4} [|\a^*_t|^4 + |\si_0(X_t)|^4]\Big] dt + g'(X_t) \si_0(X_t) dB_t = g'(X_t) \si_0(X_t) dB_t.
\eeaa
This is a $\dbP^*$-martingale. Then $\bar V_0 = Y_0 = \dbE^{\dbP^*}[Y_T] = \bar V^{\a^*}_T$. That is, $\a^*$ is an optimal control.
 
  (ii) By standard literature we also have  $V_0 = u(0,0) = g(0) = 0$. Assume by contradiction that \reff{diffusion-strong} has an optimal control $\a^*(B_\cd)$. Note that the optimal control for the Hamiltonian in \reff{diffusion-PDE} is $\sqrt{\pa^2_{xx} u(t, x)} = \si_0(x)$, then we must have  $\a^*_t (B_\cd)=\si_0(X^{\a^*}_t)$, $\dbP_0$-a.s..  Thus $X^*:= X^{\a^*}$ satisfies SDE \reff{Barlow-SDE}. 
 Since by definition $\a^*$ is $\dbF^B$-progressively measurable, we see that $X^*$ is also $\dbF^B$-progressively measurable, and hence $X^*$ is a strong solution of SDE \reff{Barlow-SDE}, contradicting with Lemma \ref{lem-Barlow}.
 \qed 

 \begin{rem}
 \label{rem2-control}
 {\rm In this example, since $\si_0$ is not differentiable in $x$, then neither is  $f$. Consequently, the stochastic maximum principle in Section \ref{sect-SMP} does not work.  
 \qed}
 \end{rem}

\subsection{Proof of Theorem \ref{thm-PDE}}
Following the arguments in Ladyzenskaja, Solonnikov \& Uralceva \cite{LSU}, we prove the theorem in four steps.

{\it Step 1.} First, for $n\ge 1$, denote
\beaa
&O_n := \{x\in \dbR^d: |x|< n\}, \q \pa O_n:= \{x\in \dbR^d: |x|=n\},&\\
&Q_n := [0, T)\times O_n,\q \pa Q_n :=( \{T\} \times O_n ) \cup ([0, T] \times \pa O_n),&\\
&g_n(t,x) := g(x)I_n(x) + [T-t] f(T, x, 0, 0)\q\mbox{and thus}\q \cL g_n(T,x) =0 ~\mbox{for}~ x\in \pa O_n,&
\eeaa
where $I_n \in C^\infty_0(\dbR^d)$ satisfying $I_n(x) = 1$ for $|x|\le n-1$ and $I_n(x) = 0$ for $|x|\ge n$. 
Next, for $k\ge 1$, define
\beaa
\si_k(t,x,y,z) := [1- I_k(z)]  I_d + I_k(z) \si(t,x, y,z),\q  f_k(t,x,y,z) := I_k(z) f(t,x, y, z).
\eeaa
Now for $k, n \ge 1$, consider the following PDE on $Q_n$:
\bea
\label{PDEkn}
\left.\ba{c}
\dis \pa_t u^k_n(t,x) + {1\over 2} \si_k^2(t,x, u^k_n, \pa_x u^k_n) : \pa^2_{xx} u^k_n + f_k(t,x, u^k_n, \pa_x u^k_n)=0,~ (t,x) \in Q_n;\\
\dis u^k_n(t,x) = g_n(t,x),\q (t,x)\in \pa Q_n.
\ea\right.
\eea
One can check that \reff{PDEkn} satisfies all the conditions in \cite{LSU} Chapter VI, Theorem 4.1, with $m=2$, $\e=0$,  $P(|z|) =0$ for $|z|\ge k$, and $\mu_1 = \mu_1(k)$ depending on $k$   in (4.6)-(4.10) there, and thus \reff{PDEkn} has a classical solution $u^k_n \in C^{1+{\b\over 2}, 2+\b}_b(Q_n \cup \pa Q_n)$ for some $\b>0$ independent of $(n, k)$.  Moreover, following the arguments of the above theorem as well as that of \cite{LSU} Chapter V, Theorem 6.1, we have
\bea
\label{uknest}
\|u^k_n \|_{C^{1+{\b\over 2}, 2+\b}_b(Q_n \cup \pa Q_n)} \le C_k,
\eea
where $C_k$ depends on $T$, $c_0, C_0, L, d$ in Assumption \ref{assam-basic}, the derivatives of the coefficients $\si, f, g$, and the index $k$, but is uniform in $n$. Now fix $k$ and send $n\to \infty$. Following the arguments of \cite{LSU} Chapter V, Theorem 8.1 and using the uniform estimate \reff{uknest}, there exists $u^k \in C_b^{1+{\b\over 2}, 2+\b}([0, T]\times \dbR^d)$ such that
\bea
\label{PDEk}
 \left.\ba{c}
\dis \pa_t u^k(t,x) + {1\over 2} \si_k^2(t,x, u^k, \pa_x u^k) : \pa^2_{xx} u^k + f_k(t,x, u^k, \pa_x u^k)=0, (t,x) \in [0, T)\times \dbR^d;\\
\dis  u^k(T,x) = g(x),\q x\in \dbR^d.
\ea\right.
\eea

{\it Step 2.} In this step we prove the first line of \reff{PDEest}. Denote
\beaa
\tilde\si_k(t,x) := \si_k\big(t,x, u^k(t,x), \pa_x u^k(t,x)\big),\q \tilde f_k(t, x, y, z) := f_k\big(t,x, y, z \tilde \si_k^{-1}(t,x)\big).
\eeaa
 By our conditions,  $\tilde \si_k$ is uniformly Lipschitz continuous in $x$, with a Lipschitz constant possibly depending on $k$, and $\tilde f_k$ is uniformly Lipschitz continuous in $(y,z)$, with Lipschitz constant uniform in $k$. By standard arguments for (strong) BSDEs, we see that
 \bea
 \label{PDErepk}
& \dis u^k(t,x) = \tilde Y^{k,t,x}_t,\q\mbox{where}&\nonumber\\
&\dis \tilde X^{k,t,x}_s = x + \int_t^s \tilde \si_k(r, \tilde X^{k,t,x}_r) dB_r,&\\
& \dis \tilde Y^{k,t,x}_s = g(\tilde X^{k,t,x}_T) + \int_s^T \tilde f_k(r, \tilde X^{k,t,x}_r, 
\tilde Y^{k,t,x}_r, \tilde Z^{k,t,x}_r) dr - \int_s^T \tilde Z^{k,t,x}_r dB_r.&\nonumber
\eea
Since $g$ and $\tilde f_k(t,x, 0,0)$ are bounded by $C_0$. It is clear that
\bea
\label{ukest}
|u^k(t,x)|\le M_0 \q\mbox{where $M_0$ depends only on $T, L$, $C_0$, $c_0$, and $d$, but not on $k$.}
\eea

We next estimate $|\pa_x u|$ under Assumption \ref{assam-basic}  (v). Note that the first case there implies $|\pa_z \si(t,x,y,z)|\le {C_0\over 1+|z|}$.  Applying \cite{LSU} Chapter VI, Theorem 3.1 on PDE  \reff{PDEkn}, with $m=2$, $\e=0$,  $P(|z|) = L$,  and $\mu_1 = C_0$ in (3.2)-(3.6) there, and passing $n\to \infty$, we obtain
\bea
\label{ukest2}
|\pa_x u^k(t,x)|\le M_1 \q\mbox{where $M_1$ depends only on $T, L$, $C_0$, $c_0$, and $d$, but not on $k$.}
\eea
In the second case that $d=1$,  denote $v^k := \pa_x u^k$. Then $v^k$ satisfies the following PDE:
\beaa
&\pa_t v^k +  {1\over 2} \tilde \si_k^2  \pa^2_{xx} v^k + \tilde b_k  \pa_x v^k + \tilde c_k v^k + \pa_x f_k (t,x,u^k, v^k),\q v^k(T,x) = \pa_x g(x),&\\
&\mbox{where}\q \tilde b_k (t,x) := \tilde\si_k \pa_x \tilde \si_k(t,x) + \pa_z f_k (t,x,u^k, v^k),\q \tilde c_k(t,x) :=  \pa_y f_k (t,x,u^k, v^k).& 
\eeaa
Note that $|\pa_x f_k|, |\pa_y f_k|, |\pa_x g|\le L$, then one may easily verify \reff{ukest2} in this case too. Now let $k \ge  M_1 + 1$, we see that $I_k(\pa_x u^k) = 1$ and thus $\f_k (t, x, u^k, \pa_x u^k) = \f(t, x, u^k, \pa_x u^k)$ for $\f = \si, f$.   That is, $u^k$ is a classical solution to the original PDE \reff{PDE}.

We finally prove the H\"{o}lder continuity of $u$ in terms of $t$.  Let $k$ be large enough and omit the subscripts $_k$ and superscripts  $^k$ in \reff{PDErepk}. Then we have the representation $u(t,x) = \tilde Y^{t,x}_t$, and $\tilde Y^{t,x}_s = u(s, \tilde X^{t,x}_s)$, $\tilde Z^{t,x}_s = \pa_x u \tilde \si (s, X^{t,x}_s)$ are bounded. For $0\le t_1 < t_2 \le T$ and $x\in \dbR^d$, by \reff{ukest2} we have
\beaa
&&|u(t_1, x) - u(t_2, x)|^2 \le C\dbE\Big[|u(t_1, x) - u(t_2, \tilde X^{t_1, x}_{t_2})|^2 + |u(t_2, \tilde X^{t_1, x}_{t_2}) - u(t_2, x)|^2\Big]\\
&&\le C\dbE\Big[|Y^{t_1,x}_{t_1} - Y^{t_1,x}_{t_2}|^2 + |\tilde X^{t_1, x}_{t_2} -\tilde X^{t_1,x}_{t_1}|^2\Big]\\
&&\le C\dbE\Big[ \int_{t_1}^{t_2} \big[|\tilde f(s, \tilde X^{t_1,x}_s,  \tilde Y^{t_1,x}_s, \tilde Z^{t_1,x}_s)|^2 + |Z^{t_1,x}_s|^2 + |\tilde \si(s, \tilde X^{t_1,x}_s)|^2\big] ds \Big] \\
&&\le C[t_2-t_1].
\eeaa
This implies the desired H\"{o}lder continuity.

{\it Step 3.}  We now prove the second line of  \reff{PDEest}. We first notice that the $C_k$ in \reff{uknest} may depend on the derivatives of the coefficients and thus \reff{uknest} does not lead to \reff{PDEest}. Instead, for any $k, n$ large, we see that $u$ satisfies the following PDE on $Q_n$ with $u$ itself as the boundary condition:
\beaa
&\pa_t u + {1\over 2} \si^2(t,x, I_k(u), I_k(\pa_x u)) : \pa^2_{xx} u + f(t,x, I_k(u), I_k(\pa_x u)) =0,\q (t,x) \in Q_n;&\\
&u(t,x) = u(t,x),\q (t,x) \in \pa Q_n.&
\eeaa
Now apply \cite{LSU} Chapter VI, Theorem 1.1, we have
\bea
\label{ukest3}
\la \pa_x u\ra^{\a}_{[0, T-\d]\times O_{n-1}} \le C_\d.
\eea
Since $n$ is arbitrary, this implies the second line of  \reff{PDEest} immediately. 

{\it Step 4.}  We finally prove (iii). First,  again by \cite{LSU} Chapter VI, Theorem 1.1, we can improve \reff{ukest3} to
\bea
\label{ukest4}
\la \pa_x u\ra^{\a}_{[0, T]\times \dbR^d} \le C_g.
\eea
Then $u$ satisfies the following linear PDE:
\bea
\label{ulinear}
\pa_t u + {1\over 2} \hat \si^2(t, x) : \pa^2_{xx} u + \hat f(t,x) =0,\q u(T,x) = g(x),
\eea
where, for $\f = \si, f$, $\hat \f(t,x) = \f(t,x, u(t,x), \pa_x u(t,x))$ is uniformly H\"{o}lder continuous. Then the estimate of $\pa^2_{xx} u$ is a classical result, see e.g. Krylov \cite{Krylov}, Theorem 8.9.2.
\qed

 \end{document}